\documentclass[english,a4paper,10pt,final]{article}
\usepackage{amsmath, amscd, amssymb, amsxtra, amsthm, array, booktabs,
  calc, color, microtype, enumitem, amsrefs, tikz-cd }
\usepackage[english]{babel} \usepackage[all]{xy}
\usepackage[colorlinks=true]{hyperref}
\begin{document}
\newcommand{\qr}[1]{\eqref{#1}} 
\newcommand{\var}[2]{\operatorname{var}_{#1}^{#2}}
\newcommand{\secref}[1]{\ref{#1}}
\newcommand{\varn}[1]{\var{#1}{}}

\newcommand{\ie}{{\it i.e.}\xspace} 
\newcommand{\eg}{{\it e.g.}\xspace}
\newcommand{\nc}{\newcommand}
\renewcommand{\frak}{\mathfrak}
\providecommand{\cal}{\mathcal}

\renewcommand{\bold}{\mathbf}
%
%
%
%
\newcommand{\gparens}[3]{{\left#1 #2 \right#3}}
\newcommand{\parens}[1]{\gparens({#1})} 
\newcommand{\brackets}[1]{\gparens\{{#1}\}} 
\newcommand{\hakparens}[1]{\gparens[{#1}]} 
\newcommand{\angleparens}[1]{\gparens\langle{#1}\rangle} 
\newcommand{\floor}[1]{\gparens\lfloor{#1}\rfloor} 
\newcommand{\ceil}[1]{\gparens\lceil{#1}\rceil}
\newcommand*{\Setof}[1]{\,\brackets{#1}\,} 
\newcommand*{\norm}[2][]{\gparens\|{#2}\|#1}
\newcommand*{\normsq}[2][]{\gparens\|{#2}\|^2#1}
%
%
\newcommand{\probopi}[3]{{#1#2}{\hakparens{\,#3\,}}}
\newcommand{\probopii}[3]{{#1#2}{\brackets{\,#3\,}}}
\newcommand{\imu}{\mu}
\newcommand{\muex}[2][\xmu]{{#1\hakparens{\,#2\,}}}
\newcommand{\mupr}[2][\xmu]{{#1\brackets{\,#2\,}}}
\renewcommand{\d}[1]{\,{d}{#1}}
\newcommand{\dmu}[2][\xmu]{\d{#1}\ifnotempty{#2}{\parens{#2}}}
\renewcommand{\Pr}[2][]{\probopii{\mathbb{P}}{#1}{#2}} 
\newcommand{\Pri}[2]{\probopii{\mathbb{P}}{_#1}{#2}} 
\newcommand{\Exp}[2][]{\probopi{\mathbb{E}}{#1}{#2}} 
\newcommand{\Exip}[2]{\probopi{\mathbb{E}}{_#1}{#2}} 
\newcommand{\EXI}[3][]{\int_{#1} {#3} \d{#2}}
\newcommand{\Zi}[2]{\probopi{Z}{#1}{#2}}
\newcommand{\probvar}[3]{{\,\mathbf{#1#2}\ifemptythenelse{#3}{}{\!\parens{#3}\,}}}
\newcommand{\Indic}[2][]{\probvar{I}{#1}{#2}} 
\newcommand{\Var}[1]{\on {Var}[{#1}]}
\newcommand{\Cov}[1]{\on {Cov}[{#1}]}
%

\numberwithin{equation}{section}

\newcommand{\theoremname}{Theorem.}
\newcommand{\corollaryname}{Corollary.}
\newcommand{\lemmaname}{Lemma.}
\newcommand{\propositionname}{Proposition.}
\newcommand{\conjecturename}{Conjecture.}
\newcommand{\definitionname}{Definition.}
\newcommand{\examplename}{Example.}
\newcommand{\remarkname}{Remark.}
\newcommand{\pfname}{Proof.}

\newenvironment{pf}{\vskip-\lastskip\vskip\medskipamount{\it\pfname}}%
                      {$\square$\vskip\medskipamount\par}

\newenvironment{pfof}[1]{\vskip-\lastskip\vskip\medskipamount{\it
    Proof of #1.}}%
                      {$\square$\vskip\medskipamount\par}

\newtheorem{thm}{Theorem}[section]
\newtheorem{theorem}[thm]{Theorem}

\newtheorem{thmnn}{Theorem}
\renewcommand{\thmnn}{}

\newtheorem{cor}[thm]{Corollary}
\newtheorem{corollary}[thm]{Corollary}
\newtheorem{cornn}{Corollary}
 \renewcommand{\thecornn}{}

\newtheorem{prop}[thm]{Proposition}
\newtheorem{proposition}[thm]{Proposition}

\newtheorem{propnn}{Proposition}
\renewcommand{\thepropnn}{}

\newtheorem{lemma}[thm]{Lemma} 
\newtheorem{lemmat}[thm]{Lemma}
\newtheorem{lemmann}{Lemma}
\renewcommand{\thelemmann}{}

\theoremstyle{definition}
\newtheorem{defn}[thm]{Definition}
\newtheorem{definition}[thm]{Definition}
\newtheorem{defnn}{Definition}
\renewcommand{\thedefnn}{}
\newtheorem{conj}[thm]{Conjecture}
\newtheorem{axiom}{Axiom}

\theoremstyle{definition}
\newtheorem{exercise}{Exercise}[subsection]
\newtheorem{exercisenn}{Exercise}
\renewcommand{\theexercisenn}{}
\newtheorem{remark}[thm]{Remark}
\newtheorem{remarks}[thm]{Remarks}
\newtheorem{remarknn}{Remark}
\renewcommand{\theremarknn}{}
\newtheorem{example}[thm]{Example} 
\newtheorem{examples}[thm]{Examples} 
\newtheorem{examplenn}{Example}[subsection]
\newtheorem{blanktheorem}{}[subsection]
\newtheorem{question}[thm]{Question}
\renewcommand{\theexamplenn}{}


\nc{\Theorem}[1]{Theorem~{#1}}
\nc{\Th}[1]{({\sl Th.}~#1)}
\nc{\Thd}[2]{({\sl Th.}~{#1} {#2})}
\nc{\Theorems}[2]{Theorems~{#1} and ~{#2}}
\nc{\Thms}[2]{({\it Thms. ~{#1} and ~{#2}})}
\nc{\Lemmas}[2]{Lemma~{#1} and ~{#2}}

\nc{\manga}[6]{({\it Thms. ~ #1, ~ #2, ~ #3,\\ ~ #4, ~ #5, ~ #6})}

\nc{\Prop}[1]{({\sl Prop.}~{#1})}
\nc{\Proposition}[1]{Proposition~{#1}}
\nc{\Propositions}[2]{Propositions~{#1} and ~{#2}}
\nc{\Props}[2]{({\sl Props.}~{#1} and ~{#2})}
\nc{\Cor}[1]{({\sl Cor.}~{#1})}
\nc{\Corollary}[1]{Corollary~{#1}}
\nc{\Corollaries}[2]{Corollaries~{#1} and ~{#2}}
\nc{\Definition}[1]{Definition~{#1}}
\nc{\Defn}[1]{({\sl Def.}~{#1})}
\nc{\Lemma}[1]{Lemma~{#1}} 
\nc{\Lem}[1]{({\sl Lem.} ~{#1})} 
\nc{\Eq}[1]{equation~({#1})}
\nc{\Equation}[1]{Equation~({#1})}
\nc{\Section}[1]{Section~{#1}}
\nc{\Sections}[1]{Sections~{#1}}
\nc{\Sec}[1]{({\sl Sec.} ~{#1})} 
\nc{\Chapter}[1]{Chapter~{#1}}
\nc{\Chapt}[1]{({\sl Ch.}~{#1})}

\nc{\Ex}[1]{{\sl Ex.}~{#1}}
\nc{\Exa}[1]{{\sl Example}~{#1}}
\nc{\Example}[1]{{\sl Example}~{#1}}
\nc{\Examples}[1]{{\sl Examples}~{#1}}
\nc{\Exercise}[1]{{\sl Exercise}~{#1}}

\nc{\Rem}[1]{({\sl Rem.}~{#1})}
\nc{\Remark}[1]{{\sl Remark}~{#1}}
\nc{\Remarks}[2]{{\sl Remarks}~{#1} and ~{#2}}
\nc{\Note}[1]{{\sl Note}~{#1}}

\nc{\Conjecture}[1]{Conjecture~{#1}}
\nc{\Claim}[1]{Claim~{#1}}

\nc \Proof{{  \it Proof. }}

\nc \Dbb{\mathbb D}
\nc{\xmu}{\mu}
\nc{\w}{\omega}
\nc{\xv}{\mbox{\boldmath$x$}}
\nc{\uv}{\mbox{\boldmath$u$}}
\nc{\xiv}{\mbox{\boldmath$\xi$}}
\nc{\bbeta}{\mbox{\boldmath$\beta$}}
\nc{\balpha}{\mbox{\boldmath$\alpha$}}
\nc{\bgamma}{\mbox{\boldmath$\gamma$}}
\nc{\bdelta}{\mbox{\boldmath$\delta$}}
\nc{\bepsilon}{\mbox{\boldmath$\epsilon$}}

\newcommand{\ZZ}{{\mathbb Z}}
\newcommand{\RR}{{\mathbb R}} 

\nc \Ab{{\ensuremath{\bold A}}}
\nc \ab{{\ensuremath{\bold a}}}
\nc \bb{{\ensuremath{\bold b}}}
\nc \cb{{\ensuremath{\bold c}}}
\nc \db{{\ensuremath{\bold d}}}
\nc \Bb{{\ensuremath{\bold B}}}
\nc \Gb{{\ensuremath{\bold G}}}
\nc \Qb{{\ensuremath{\bold Q}}}
\nc \Rb{{\ensuremath{\bold R}}} \nc \Cb{{\ensuremath{\bold C}}} 
\nc \Eb{{\ensuremath{\bold E}}}
\nc \eb{{\ensuremath{\bold e}}}
\nc \Db{{\ensuremath{\bold D}}}
\nc \Fb{{\ensuremath{\bold F}}}
\nc \ib{{\ensuremath{\bold i}}}
\nc \jb{{\ensuremath{\bold j}}}
\nc \kb{{\ensuremath{\bold k}}}
\nc \Kb{{\ensuremath{\bold K}}}
\nc \nb{{\ensuremath{\bold n}}}
\nc \rb{{\ensuremath{\bold r}}}
\nc \Ob{{\ensuremath{\bold O}}}
\nc \Pb{{\ensuremath{\bold P}}}
\nc \pb{{\ensuremath{\bold p}}}
\nc \qb{{\ensuremath{\bold q}}}
\nc \SPb{{\ensuremath{\bold {SP}}}}
\nc \Zb{{\ensuremath{\bold Z}}} 
\nc \zb{{\ensuremath{\bold z}}} 
\nc \gb{{\ensuremath{\bold g}}} 
\nc \fb{{\ensuremath{\bold f}}} 
\nc \ub{{\ensuremath{\bold u}}} 
\nc \vb{{\ensuremath{\bold v}}} 
\nc \yb{{\ensuremath{\bold y}}} 
\nc \xb{{\ensuremath{\bold x}}} 
\nc \tb{{\ensuremath{\bold t}}} 
\nc \Xb{{\ensuremath{\bold X}}} 
\nc \Yb{{\ensuremath{\bold Y}}} 
\nc \xib{{\ensuremath{\bold \xi}}} 
\nc \Nb{{\ensuremath{\bold N}}} 
\nc \Hb{{\ensuremath{\bold H}}} 
\nc \wb{{\ensuremath{\bold w}}} 
\nc \Wb{{\ensuremath{\bold W}}} 
\nc \syz{{\mathbf {syz}}}
\nc \bnoll{{\ensuremath{\bold 0}}} 

\nc \mf{\frak m} 
\nc \mh{\hat{\mf}} 
\nc \nf{\frak n}
\nc \Of{\frak O}
\nc \of{\frak o}
\nc \rf{\frak r}
\nc \tf{\frak t}
\nc \mufr{{\mathbf \mu}}
\nc \hf{\frak h} 
\nc \qf{\frak q} 
\nc \bfr{\frak b} 
\nc \kfr{\frak k} 
\nc \pfr{\frak p} 
\nc \af{\frak a }
\nc \cf{\frak c }
\nc \sfr{\frak s} 
\nc \ufr{\frak u} 
\nc \g{\frak g} 
\nc \gA{\g_{\Ao}} 
\nc \lfr{\frak l}
\nc \afr{\frak a}
\nc \gfh{\hat {\frak g}}
\nc \gl{\frak { gl }}
\nc \Sl{\frak {sl}}
\nc \SU{\frak {SU}}
\nc{\Homf}{\frak{Hom}}

\nc{\Av}{\overline{A}}

\newcommand{\on}{\operatorname}
\newcommand{\adj}[2]{ \overset{#1}{\underset{#2}\rightleftarrows }}
\nc\hankel{\on {Hankel}}
\nc\row{\on {row\ }}
\nc\nullity{\on {nullity }}
\nc\col{\on{col}}
\nc\rowm{\on {Row \ }}
\nc\loc{\on {lc \ }}
\nc\lcm{\on {lcm \ }}
\nc\nullo{\on {null\ }}
\nc\Nul{\on {Nul}}
\nc \Ann {\on {Ann }}
\nc \Ass {\on {Ass}}
\nc \Coker {\on {Coker}}
\nc \Co{\on C}
\nc \Homo{\on {Hom}}
\nc \Galo{\on {Gal}}
\nc \Ker {\on {Ker}}
\nc \Lef {\on {Lef}}
\nc \Leff {\on {Leff}}
\nc \omod{\on{mod}}
\nc \No {\on N}
\nc \NN {\on {NN}}
\nc \NGo {\on {NG}}
\nc \Oo {\on O}
\nc \ch {\on {ch}}
\nc \cd {\on {cd}}
\nc \rko {\on {rk}}
\nc \Sing {\on {Sing\ }}
\nc \Reg {\on {Reg}}
\nc \RG{\on {R\Gamma}}
\nc \CoI {\on {CI}}
\nc \CoM {\on {CM}}
\nc \Gor {\on {Gor}}
\nc \Type {\on {Type}}
\nc \can {\on {can}}
\nc \Top {\on {T}}
\nc \Tr {\on {Tr}}
\nc \Norm {\on {Norm}}
\nc \rel {\on {rel}}
\nc \tr {\on {tr}}
\nc \sgn {\on {sgn }}
\nc \trdeg {\on {tr.deg}}
\nc \codim {\on {codim }}
\nc \coht {\on {coht}}
\nc \divo {\on {div \ }}
\nc \rot {\on {rot }}
\nc \coh {\on {coh}}
\nc \Clo {\on {Cl}}
\nc \Divo {\on {Div}}
\nc \embdim{\on {embdim}}
\nc \ord{\on {ord}}
\nc \ed{\on {ed}}
\nc \embcodim{\on {embcodim  }}
\nc \qcoh {\on {qcoh}}
\nc \grad {\on {grad}}
\nc \grade {\on {grade}}
\nc \hto {\on {ht}}
\nc \depth {\on {depth}}
\nc \prof {\on {prof}}
\nc \reso{\on {res}}
\nc \Reso{\on {Res}}
\nc \ind{\on {ind}}
\nc \prodo{\on {prod}}
\nc \coind{\on {coind}}
\nc \Con{\on {Con}}
\nc \Crit{\on {Crit}}
\nc \Der{\on {Der}}
\nc \Des{\on {Des}}
\nc \Char{\on {Char}}
\nc \Ch{\on {Ch}}

\nc \Ext{\on {Ext}}
\nc \Eo{\on {E}}
\nc \End{\on {End}}
\nc \ad{\on {ad}}
\nc \Ad{\on {Ad}}
\nc \gr{\on {gr}}
\nc \Fo{\on {F}}
\nc \Gr{\on {Gr}}
\nc \Go{\on {G}}
\nc \GFo{\on {GF}}
\nc \Glo{\on {Gl}}
\nc \PGlo{\on {PGl}}
\nc \Ho{\on {H}}
\nc \CMo{\on {\CM}}
\nc \SCM{\on {SCM}}
\nc \rig{\on {right}}
\nc \lef{\on {left}}
\nc \hol{\on {hol}}
\nc{\sgd}{\on{sgd}}
\nc \supp{\on {supp}}
\nc \ssupp{\on {s-supp}}
\nc \singsupp{\on {singsupp}}
\nc \msupp{\on {msupp}}
\nc \spec{\on {spec}}
\nc \spano{\on {span }}
\nc \Span{\on {Span }}
\nc \Max{\on {Max}}
\nc \Mat{\on {Mat}}
\nc \Min{\on {Min}}
\nc \nil{\on {nil}}
\nc \Mod{\on{Mod}}
\nc \Rad {\on {Rad}}
\nc \rad {\on {rad}}
\nc \rank {\on {rank}}
\nc \range {\on {range}}
\nc \Slo{\on {SL}}
\nc \soc {\on {soc}}
\nc \dt {\on {dt}_Z}
\nc \Irr {\on {Irr}}
\nc \Reo {\on {Re}}
\nc \Imo {\on {Im}}
\nc \SSo{\on {SS}}
\nc \lub{\on {lub}}
\nc \gldim{\on {gl.d.}}
\nc \length{\on {length}}
\nc \pdo{\on {p.d.}} 
\nc \cdo{\on {cd}} 
\nc \fdo{\on {f.d.}} 
\nc \ido{\on {i.d.}} 
\nc \dSSo{\dot {\SSo}}
\nc \So{\on S}
\nc \SOo{\on{ SO}}
\nc \Io{\on I}
\nc \Jo{\on J}
\nc \jo{\on j}
\nc \Ko{\on K}
\nc \PBW{\Ac_{PBW}}
\nc \Ro{\on R}
\nc \To{\on T}
\nc \Ao{\on A}

\nc \Do{{\on D}}
\nc \Bo{\on B}
\nc \Po{\on P}
\nc \Qo{\on Q}
\nc \Zo{\on Z}
\nc \Uo{\on U}
\nc \wt{\on {wt}}
\nc \Uoh{\hat {\Uo}}
\nc \Lo{\on L}
\nc \Loc{\on {Loc}}
\nc{\dop}{\on d}
\nc{\eo}{\on e}
\nc{\ado}{\on{ad}}
\nc{\Tot}{\on{Tot}}
\nc{\Aut}{\on{Aut}}
\nc{\sinc}{\on {sinc}}

%
%
\nc{\overrightleftarrows}[2]{\overset{#1}{\underset{#2}{\rightleftarrows}}}

\nc{\CCF}{\cal{CF}}
\nc{\CDF}{\cal{DF}}
\nc{\CHC}{\check{\cal C}}

\nc{\Cone}{\on{Cone}}
\nc{\dec}{\on{dec}}
\nc{\Diff}{\on{Diff}}
\nc{\dirlim}{\underset{\to}{\on{lim}}}
\nc{\dpar}{\partial}
\nc{\dlog}{\on{dlog}}
\nc{\GL}{\on{GL}}
\nc{\glo}{\on{gl}}
\nc{\CGr}{\cal{G}r}
\nc{\pr}{\on{pr}}
\nc{\semid}{|\!\!\!\times}
\nc{\Hom}{\on{Hom}}
\nc \RHom{\on {RHom}}

\nc \Proj{\mathrm {Proj\ }}
\nc \proj{\mathrm {proj}}
\nc{\Id}{\on{Id}}
\nc{\id}{\on{id}}
\nc{\Ima}{\on{Im}}
\nc{\invtimes}{\underset{\gets}{\otimes}}
\nc{\invlim}{\underset{\gets}{\on{lim}}}
\nc{\Lie}{\on{Lie}}
\nc{\re}{\on{Re }}
\nc{\Pic}{\on{Pic }}
\nc{\LPic}{\on{LPic }}
\nc{\Sch}{\on{Sch}}
\nc{\Sh}{\on{Sh}}
\nc{\Set}{\on{Set}}
\nc{\spo}{\on{sp\  }}
\nc{\Spec}{\on{Spec}}
\nc{\mSpec}{\on{mSpec}}
\nc{\Specb}{\bold {Spec}\ }
\nc{\Projb}{\bold {Proj}}
\nc{\Specan}{\on{Specan}}
\nc{\Spo}{\on{Sp}}
\nc{\Mpo}{\on{Mp}}
\nc{\Spf}{\on{Spf}}
\nc{\sym}{\on{sym}}
\nc{\symm}{\on{symm}}
\nc{\rop}{\on{r}}
\nc{\Td}{\on{Td}}
\nc{\Tor}{\on{Tor}}


\nc{\Alg}{\on {Alg}}
\nc{\Artin}{\cal{A}rtin}
\nc{\Dgcoalg}{\cal{D}gcoalg} \nc{\Dglie}{\cal{D}glie}
\nc{\Ens}{\cal{E}ns} \nc{\Fsch}{\cal{F}sch}
\nc{\Groupoids}{\cal{G}roupoids}
\nc{\Holie}{\cal{H}olie}
\nc{\Mor}{\cal{M}or}

\nc \Dd{\mathbb D}
\nc{\CF}{\ensuremath{\cal{F}}}
\nc \Kc{{\ensuremath{\cal K}}}
\nc{\Kzind}[4]{{\ensuremath{{\mathcal K^{#4e} (#3)}_{#1 \rightarrow
        #2}}}}
\nc{\Kz}[3]{{\ensuremath{{\mathcal K^\bullet (#3)}_{#1 \rightarrow
        #2}}}}
\nc{\Kzd}[2]{{\ensuremath{{\mathcal K}^\bullet_{#1 \rightarrow #2}}}}
\nc \Lc{{\ensuremath{\cal L}}}
\nc \lcc{{\mathcal l}} 
\nc \CC{{\ensuremath{\cal C}}} 
\nc \Cc{{\ensuremath {\cal C}}}
\nc{\Dl}[2]{{\ensuremath{{\mathcal D}_{#1 \leftarrow #2}}}}
\nc{\Dr}[2]{{\ensuremath{{\mathcal D}_{#1 \rightarrow #2}}}}
\nc \Dc{\ensuremath{\mathcal D}}
\nc \Ac{{\ensuremath{\cal A}}} 
\nc \Bc{{\ensuremath{\cal B}}}
\nc \Ec{{\ensuremath{\cal E}}}
\nc \Fc{{\ensuremath{\cal F}}}
\nc \Mcc{{\ensuremath{\cal M}}} 
\nc \hM{\hat{\Mcc}} 
\nc \bM{\bar {\Mcc}} 
\nc\hbM{\hat{\bar \Mcc}}  
\nc \Nc{{\ensuremath{\cal N}}}
\nc \Hc{{\ensuremath{\cal H}}} 
\nc \Ic{{\ensuremath{\cal I}}} 
\nc \Jc{{\ensuremath{\cal J}}} 
\nc \Oc{\ensuremath{{\cal O}}}
\nc \Och{\hat{\cal O}} 
\nc \Sc{{\ensuremath{{\cal S}}}}
\nc \Tc{\ensuremath{{\cal T}}} 
\nc \Vc{{\ensuremath{{\cal V}}}} 
\nc{\CA}{{\ensuremath{{\cal A}}}}
\nc{\CB}{{\ensuremath{{\cal B}}}}
\nc{\Fcc}{\ensuremath{{\cal F}}}
\nc{\Gc}{{\ensuremath{{\cal G}}}}
\nc{\CH}{\ensuremath{\mathcal H}}
\nc{\CI}{{\ensuremath{{\cal I}}}}
\nc{\CM}{{\ensuremath{{\cal M}}}}
\nc{\CN}{{\ensuremath{{\cal N}}}}
\nc{\CO}{{\ensuremath{{\cal O}}}}
\nc{\Rc}{{\ensuremath{{\cal R}}}}
\nc{\Pc}{{\ensuremath{{\cal P}}}}
\nc{\CT}{{\ensuremath{\mathcal T}}}
\nc{\CU}{\ensuremath{{\cal U}}}
\nc{\Uc}{\ensuremath{{\cal U}}}
\nc{\Yc}{\ensuremath{{\cal Y}}}
\nc{\CV}{\ensuremath{{\cal V}}}
\nc{\CZ}{\ensuremath{{\cal Z}}}
\nc{\Homc}{\ensuremath{{\cal {Hom}}}}

\nc{\fa}{\frak{a}}
\nc{\fA}{\frak{A}}
\nc{\fg}{\frak{g}}
\nc{\fh}{\frak{h}}
\nc{\fI}{\frak{I}}
\nc{\fK}{\frak{K}}
\nc{\fm}{\frak{m}}
\nc{\fP}{\frak{P}}
\nc{\fS}{\frak{S}}
\nc{\ft}{\frak{t}}
\nc{\fX}{\frak{X}}
\nc{\fY}{\frak{Y}}


\nc{\bF}{\bar{F}}
\nc{\bCP}{\bar{\cal{P}}}
\nc{\bmbox}{\mbox{\bf{m}}}
\nc{\bT}{\mbox{\bf{T}}}
\nc{\hB}{\hat{B}}
\nc{\hC}{\hat{C}}
\nc{\hP}{\hat{P}}
\nc{\htest}{\hat P}


\nc{\nen}{\newenvironment}
\nc{\ol}{\overline}
\nc{\unl}{\underline}
\nc{\ra}{\to}
\nc{\lla}{\longleftarrow}
\nc{\lra}{\longrightarrow}
\nc{\Lra}{\Longrightarrow}
\nc{\Lla}{\Longleftarrow}
\nc{\Llra}{\Longleftrightarrow}
\nc{\hra}{\hookrightarrow}
\nc{\iso}{\overset{\sim}{\lra}}

\nc{\dsize}{\displaystyle}
\nc{\sst}{\scriptstyle}
\nc{\tsize}{\textstyle}
%

%
%
%
\nen{exa}[1]{\label{#1}{\bf Example.\ } }{}


\nen{rem}[1]{\label{#1}{\em Remark.\ } }{}
\nen{exer}[1]{\label{#1}{\em Exercise.\ } }{}



\newtheorem{sats}{Sats}[section]
\newtheorem{satsnn}{Sats}
\renewcommand{\thesatsnn}{}

\newtheorem{kor}[sats]{Korollarium}
\newtheorem{korollarium}[sats]{Korollarium}
\newtheorem{korollariumnn}{Korollarium}
\newtheorem{kornn}{Korollarium}
\renewcommand{\thekornn}{}

\newtheorem{sprop}[sats]{Proposition}
\newtheorem{sproposition}[sats]{Proposition}

\newtheorem{spropnn}{Proposition}
\renewcommand{\thepropnn}{}

\newtheorem{slemma}[sats]{Lemma} 
\newtheorem{slemmat}[sats]{Lemma}
\newtheorem{slemmann}{Lemma}
\renewcommand{\theslemmann}{}

\theoremstyle{definition}
\newtheorem{sdefn}[sats]{Definition} 
\newtheorem{sdefinition}[sats]{Definition} 
\newtheorem{formodan}[sats]{Förmodan}

\theoremstyle{remark}
\newtheorem{ovning}[sats]{Övning}
\newtheorem{ovningnn}{Övning}
\newtheorem{svarnn}{Svar}
\renewcommand{\thesvarnn}{}
\newtheorem{ovningarnn}{Övningar}
\renewcommand{\theovningarnn}{}
\renewcommand{\theovningnn}{}
\newtheorem{anm}[sats]{Anmärkning}
\newtheorem{anmnn}{Anmärkning}
\renewcommand{\theanmnn}{}
\newtheorem{exempel}[sats]{Exempel} 
\newtheorem{exempelnn}{Exempel} 
\renewcommand{\theexempelnn}{}


\nc{\Sats}[1]{Sats~\ref{#1}}
\nc{\Sa}[1]{({\sl Sa.}~\ref{#1})}
\nc{\Kor}[1]{({\sl Kor.}~\ref{#1})}
\nc{\Korollarium}[1]{Korollarium~\ref{#1}}
\nc{\Exe}[1]{{\sl Exempel}~\ref{#1}}
\nc{\Anm}[1]{{\sl Anmärkning}~\ref{#1}}
\nc{\Not}[1]{{\sl Not}~\ref{#1}}

\nc{\Formodan}[1]{Förmodan~\ref{#1}}
\nc{\Pastaende}[1]{Påstående~\ref{#1}}

\newenvironment{bv}{\vskip-\lastskip\vskip\medskipamount{\it Bevis.}}%
                      {$\square$\vskip\medskipamount\par}


\title{The Zariski-Lipman conjecture for complete intersections}
\author{Rolf K\"allstr\"om}
\maketitle

\footnotetext[1]{2000 Mathematics Subject Classification: {Primary:
    14A10, 32C38; Secondary: 17B99 (Secondary)}} 
\begin{abstract} The tangential ramification locus
  $B_{X/Y}^t\subset B_{X/Y}$ is the subset of points in the
  ramification locus where the sheaf of relative vector fields
  $T_{X/Y}$ fails to be locally free. It was conjectured by Zariski
  and Lipman that if $V/k$ is a variety over a field $k$ of
  characteristic $0$ and $B^t_{V/k}= \emptyset$, then $V/k$ is smooth
  (=regular). We prove this conjecture when $V/k$ is a locally
  complete intersection. We prove also that $B_{V/k}^t= \emptyset$
  implies $\codim_V B_{V/k}\leq 1 $ in positive characteristic, when
  $V/k $ is the fibre of a flat morphism satisfying generic
  smoothness.
\end{abstract}

\section{Introduction}
Let $\pi: X\to Y$ be a morphism of n{\oe}therian schemes which is
locally of finite type, $\Omega_{X/Y}$ its sheaf of K{\"a}hler
differentials, and $T_{X/Y}= Hom_{\Oc_X}(\Omega_{X/Y}, \Oc_X)$ the
sheaf of relative tangent vector fields. We have the inclusion of the
tangential ramification locus in the ramification locus
\begin{displaymath} B^t_\pi= \{x\in X \ \vert \ T_{X/Y,x} \text{ is
    not free}\} \subset B_\pi =B_{X/Y}= \{x \in X \ \vert \
  \Omega_{X/Y,x} \text{ is not free}\}.
\end{displaymath}
Define as in \cite[Definitions 17.1.1 and 17.3.1]{EGA4:4} a morphism
$\pi$ to be formally smooth at a point $x$ in $X$ if the induced map
of local rings $\Oc_{Y,\pi(x)}\to \Oc_{X,x}$ is formally smooth, and
that $\pi$ is smooth at $x$ if it is locally finitely presented and
formally smooth; say also that $\pi$ is smooth if it is smooth at all
points in $X$. In the light of the fact that the Jacobian criterion,
namely that $B_\pi=\emptyset$, goes a long way to implying that the
morphism $\pi$ is smooth \Thms{\ref{smooth-flat}} {\ref{th-smooth}},
it is a natural to ask, with Zariski and Lipman \cite{lipman:freeder},
what are the implications of $B^t_{\pi}=\emptyset$? The example
$X=\Spec A[x]/(x^2) \to Y=\Spec A$, i.e. the scheme of dual numbers
over a commutative ring $A$, shows that if we want $\pi$ to be smooth,
the condition $B^t_\pi= \emptyset$ needs at least to be supplemented
with the condition that the rank of $T_{X/Y}$ equals the relative
dimension at each point in $X$, which can be imposed by assuming that
$X/Y$ is smooth at generic points in $X$. It is a remarkable fact that
although $T_{X/Y}$ cannot even directly detect torsion in
$\Omega_{X/Y}$, it turns out that these conditions combined imply
$B_\pi=\emptyset$ (and hence imply that $\pi$ is smooth) in
interesting cases in characteristic $0$. Already the result that
$B_{V/k}^t=0$ implies smoothness when $V/k$ is a curve over a field of
characteristic $0$, due to Lipman [loc. cit], is surprising and
non-trivial (see \Proposition{\ref{lipman-normal}(2)}). In positive
characteristic it is easy to see that smoothness at points of height
$\leq 1$ does not follow from $B^t_\pi =\emptyset$, so one could
perhaps add the assumption $\codim_X B_\pi \geq 2$, but this is still
not enough. What is needed is a condition on the discriminant locus
$D_\pi = \pi(B_\pi)$.

The last two sections are not used for the main result. One is a
discussion of the Jacobian criterion of smoothness for morphisms that
are locally of finite type, with the intention to amend/explain a
passage in \cite{EGA4:4}, and the other is about expressing conditions
under which a non-smooth map is submersive; it is also a kind of
converse to a result due to Lichtenbaum and Schlessinger
\Prop{\ref{lipman-normal}}(1).

Before the main results are presented we
describe some terminology.

\medskip {\it Generalities}: All schemes are assumed to be
n{\oe}therian and we use the notation in EGA, but see also \cite[\S
5]{matsumura} and \cite{hartshorne}. The height $\hto_X (x)$ of a
point $x$ in $X$ is the same as the Krull dimension of the local ring
$\Oc_{X,x}$ at $x$, and the dimension of $X$ is defined as
$\dim X = \sup \{\ \hto (x)\ \vert \ x\in X\}$. There is a partial
order on $X$ where a point $x$ is greater $x'$ if $x' \in \{x\}^-$,
where $\{x\}^-$ is the closure of $x$; in other words, $x$ specialises
to $x'$ and we write $x \rightsquigarrow x'$ (see \cite[p.
93]{hartshorne}); in particular $\hto_X(x')\geq \hto_X(x)$. A point
$x$ in a subset $T$ of $X$ is {\it maximal} if it is maximal for this
partial order, that is, if $x' \in T$ and $x' \rightsquigarrow x$,
then $x'=x$. Denote by $\Max (T)$ the set of maximal points of $T$, so
$\Max(X)$ consists of points of height $0$. A property on $X$ is {\it
  generic} if it holds for all points in $\Max (X)$.
Put \begin{eqnarray*}
      \codim^+_X T &=& \sup \{\ \hto_X (x) \  \vert \ x \in \Max (T) \},\\
      \codim^-_X T &=& \inf \{\ \hto_X (x) \ \vert \ x \in \Max (T)\},
 \end{eqnarray*}
 so $\codim^-_{X} T \leq \hto (x)\leq \codim^+_X T$ when $x\in \Max
 (T)$. If $T$ is the empty set, put $\codim_X^{+}T=-1$ and $\codim_X^-
 T=\infty$, since we are interested in lower and higher bounds on
 $\codim_X ^\pm T$, respectively.  For a coherent $\Oc_X$-module $M$,
 the stalk at a point $x$ is denoted $M_x$ and we put $\depth_T M =
 \inf \{\depth M_x \ \vert \ x\in T\}$, where $\depth M_x$ is the
 maximal length of a regular sequence in the ideal $\mf_{X,x}\subset
 \Oc_{X,x }$; we also put $\depth (x) = \depth (\Oc_{X,x})$. The fibre
 $X_y$ over a point $y$ in $Y$ is the fibre product $\Spec
 k_{Y,y}\times_Y Y$, where $k_{Y,y}$ is the residue field at $y$. Let
 $\hto_X(x/x')$ be the maximal length $n$ of a chain of distinct
 points $x' \rightsquigarrow  x_1 \rightsquigarrow  \cdots \rightsquigarrow  x_n = x$, where $x_i$ specialises
 to $x_j$ when $j  >i$. The dimension at a point $x$ is the maximal
 length of such chains that pass through $x$:
\begin{displaymath}
  \dim_x X = \sup \{ \hto_X (x_1/x) \ \vert \  x
  \rightsquigarrow x_1 \} + \hto_X(x).
\end{displaymath}
This number equals
$ \sup \{ \ \hto_X (x_1) \ \vert \ x_1 \in X, \quad x \rightsquigarrow x_1
\}$ when $X$ is equidimensional and catenary at all specialisations of
$x$. See \citelist{\cite{EGA4:1}*{Definition
    14.2.1}\cite{EGA4:2}*{Proposition 5.1.4}\cite{matsumura}*{Lemma 2,
    p. 250}} and \Theorem{\ref{hartshorne}}. A minimal point $x $ in
$T$ is one such that if $x$ specialises to another point $x'$ in $T$,
then $x= x'$. Then if $T$ is a closed subset it follows that $\{x\}$
is a closed subset of $X$, since $X$ is n{\oe}therian and hence
quasi-compact; therefore minimal points are the same as closed points
in a closed subset. We define the {\it relative dimension} $d_{X/Y,x}$
of $\pi$ at a point $x\in X$ as the infimum of the dimension of the
vector space of K{\"a}hler differentials at all maximal points $\xi$
that specialise to $x$, i.e.
 \begin{displaymath}
   d_{X/Y,x}=\inf \{
   \dim_{k_{X,\xi}}k_{X,\xi}\otimes_{\Oc_{X,\xi}}\Omega_{X/Y,\xi} \
   \vert \ \xi \rightsquigarrow x  , \quad  \xi \in \Max (X)\}.
 \end{displaymath}
 To understand this number it is useful recall that
\begin{align*}
\dim_{k_{X,\xi}}k_{X,\xi}\otimes_{\Oc_X,\xi
 }\Omega_{X/Y,\xi} &= \dim_{k_{X_{\pi(\xi)}, \xi}}
 \Omega_{X_{\pi(\xi)}/k_{Y, \pi(\xi)}}\\ 
&= \dim_{k_{X, \xi}}
 k_{X,\xi}\otimes_{\Oc_{X,\xi}}\Omega_{\Oc_{X, \xi}/\Oc_{Y,\pi(\xi)}};
\end{align*}
see \Proposition{\ref{base-change-differentials}} for the first
equality, but note that in general the numbers $d_{X/Y,x}$ and $\dim_x
X_{\pi(x)} $ are not equal. On the other hand, if $\pi$ is flat at
$x$, then $\hto_ {X_{\pi(x)}}(x) = \dim \Oc_{X,x} - \dim
\Oc_{Y,\pi(x)}$, which equals $\dim_x X_{\pi(x)}$ when $x$ is a closed
point in the fibre $X_{\pi(x)}$, and if morover $\pi$ is smooth at all
points $\xi\in \Max (X)$ that specialise to $x$, then $d_{X/Y,x}=
\dim_x X_{\pi(x)}$.

When we say that $A\to B$ is a homomorphism of rings, we mean a local
homomorphism of local n{\oe}therian rings. We let $\mf_A$ and $k_A$
denote the maximal ideal and residue field of $A$, respectively.

Say that a dominant morphism $X\to Y$ locally of finite type, is a
differentially complete intersection (d.c.i.) if the projective
dimension $\pdo \Omega_{X/Y,x}\leq 1$ at each point $x$. Ferrand and
Vasconcelos, has proven that $X/Y$ is a d.c.i. if $X/Y$ is generically
smooth (in $X$) and if it locally can be factorised
$X \to \bar X \to Y$, where $X \to \bar X$ is a regular immersion and
$\bar X/Y$ is smooth. For example, if $X/S$ and $Y/S$ are smooth and
$X/Y$ is generically smooth, then $X/Y$ is a d.c.i. See
\cite{kallstrom:branchpurity} for a discussion of this notion.

\begin{theorem}\label{Z-L} Let $X/S$ and $Y/S$ be
  n{\oe}therian $S$-schemes which are locally of finite type, $X$ is
  Cohen-Macaulay, $B_{X/S}= \emptyset$, and $\pi:X/S \to Y/S$ be a
  flat d.c.i. $S$-morphism, and $B_{X/S}= \emptyset$ (e.g. $X/S$ is
  smooth). Assume that $\codim_Y^- D_{\pi}\geq 1$.
  \begin{enumerate}
    \item If $x\notin B_{X/Y}^t $, and $\xi \rightsquigarrow x$, where
      $\xi \in \Max (B_{X_y/k_{Y,y}})$, then $\hto_{X_y}(\xi )\leq 1$.
    \item Let $y$ be a point in $Y$ such that $\Oc_{Y,y}$ is regular.
      If $\codim^-_{X_y} B_{X_y/ k_Y,y} \geq 1$, then
      $B^t_{X/Y}\cap X_y = B^t_{X_y/k_{Y,y}}$.
  \end{enumerate}

\end{theorem}
\begin{remark}\label{gensmooth}
  The condition of generic smoothness on $Y$, i.e.
  $\codim_Y^- D_{X/Y}\geq 1$, is satisfied when $\Oc_{X,x}$ is regular
  and the extension of residue fields $k_{X,x}/k_{Y,y}$ is separable
  for all closed points $x\in X_y$, when $ y \in \Max (Y)$. A
  \hyperlink{assertion}{proof} of this assertion is included at the
  end.
\end{remark}
\begin{corollary}\label{Z-L-cor}
  Let $V/k$ be a variety defined by a regular sequence $\{f_1, \dots ,
  f_r\} $ in some polynomial ring $k[X_1, \cdots , X_n]$ and assume
  that $T_{V/k}$ is locally free.
  \begin{enumerate}\item 
    If $\Char k =0$, then $V/k$ is smooth.
  \item If $\Char k >0$, assume moreover that the ring $k(f_1, \dots ,
    f_r)\otimes_{k[f_1, \dots , f_r]} k[X_1, \dots , X_n]$ is smooth
    over the field $k(f_1, \dots , f_r)$.  Then
    \begin{displaymath}
\codim^+_VB_{V/k}    \leq 1.
  \end{displaymath}
 \end{enumerate}

\end{corollary}
\begin{remarks}
  \begin{enumerate}
  \item The assumption in \Corollary{\ref{Z-L-cor}} that $V/k$ is
    defined by a regular sequence is used to infer that a morphism
    $A\to B$ in the proof is flat (the local flatness criterion).
    Conversely, if $A\to B$ is flat then the fibre $V$ is a complete
    intersection \cite{avramov:flatlci} (see also
    \cite{kallstrom:branchpurity}).
  \item Zariski and Lipman \cite{lipman:freeder} stated their
    conjecture only for varieties over fields of characteristic $0$.
    We can ``explain'' the positive characteristic counterexample in
    [loc. cit, \S 7,b)]. The surface $V= V( XY-Z^p)\subset \Ab_k^3$
    over a perfect field $k$ of characteristic $p>0$ is normal and
    $T_{V/k}$ is locally free. By normality and since $k$ is perfect,
    $V$ is smooth at all points of height $\leq 1$, in accordance with
    \Corollary{\ref{Z-L-cor}}. Since $V/k$ is not smooth at the
    origin, \Theorem{\ref{Z-L}} implies that if $V$ is the fibre of a
    flat family of surfaces $X\to Y$, where $X/k$ and $Y/k$ are
    smooth, then $X/Y$ cannot be generically smooth in $Y$. For
    example, the hypersurface $X=V(t-XY-Z^p)\subset \Ab^4_k$ is smooth
    over $k$, the morphism $\pi: X\to Y =\Ab^1_k $ induced by the
    projection to the $t$-coordinate is flat, and $T_{X/Y}$ is locally
    free. However, $\pi$ is not generically smooth on $Y$ since the
    field extension $k_{X,x}/k_{Y,\pi(x)}$ is not separable when $x$
    is the maximal point in $X$.
 \end{enumerate}
\end{remarks}
Scheja and Storch \cite{scheja-storch:eigenschaften} proved
\Corollary{\ref{Z-L-cor}}(1) when $V/k$ is a hypersurface in
characteristic $0$ (their proof is included); Moen
\cite{moen:free_derivation} proved it when $V$ is a homogeneous
complete intersection; Hochster \cite{hochster:zariski-lipman} proved
it when $V$ is the spectrum of a graded ring, and he even attempted to
find a counter-example when $V/k$ is a locally complete intersection
surface. Platte \cite{platte:zl} found an elementary proof in the
graded case, applicable also to analytic algebras. Using the existence
of weakly submersive resolutions of singularities (work of Hironaka)
combined with a detailed study of mixed Hodge structures, Straten and
Steenbrink \cite{straten-steenbrink} argue that
$\codim^+_X B_{X/\Cb}\leq 2$ when $T_{X/\Cb}$ is locally free and $X$
is an analytic space of dimension at least 3 with at most isolated
singularities; this was extended to non-isolated singularities by
Flenner \cite{flenner:extendability} assuming
$\codim^+_X B_{X/\Cb}\leq 2$. This can be compared to Lichtenbaum and
Schlessinger's result \cite [Prop. 5.2]{lipman:freeder}, described in
\Proposition{\ref{lipman-normal}}(1).

\begin{remark}
  First I want to thank a referee for helping me make
  \cite{kallstrom:zl} more accessible. However, the present paper is a
  rather complete rewrite of [loc cit], whose main problem is that it
  contains a serious mistake (see [loc cit, \Proposition{2.2},(ii) in
  the proof]), and its correction required new ideas. I want to thank
  very much Mathias Schulze for posing good questions leading to the
  discovery of the mistake.
\end{remark}

\section{Base change for relative tangent vector fields}\label{basechange}
Let $\pi: X\to Y$ be a morphism of connected n{\oe}therian schemes
which is locally of finite type and generically smooth, so
$\Omega_{X/Y,x}$ is free of rank $d_{X/Y,x}$ when $x$ is a maximal
point. Put $d_{X/Y}= \inf_{x\in X} d_{X/Y,x}$.  The ramification scheme
$B^{(i)}_{X/Y}$, $i=0, \dots , $ is defined by the Fitting ideal
$F_{d_{X/Y}+i}(\Omega_{X/Y})$, and we put
$B_\pi=B_{X/Y}=B^{(0)}_{X/Y}= V(F_{d_{X/Y}}(\Omega_{X/Y}))$, which is
the locus of points where $\Omega_{X/Y}$ is not free of rank
$d_{X/Y}$.  Similarly, $B^{t,(i)}_{X/Y}= V(F_{d_{X/Y}+i}(T_{X/Y}) )$,
and the tangential ramficiation scheme - the locus of points where $T_{X/Y}$
is not free - is $B^t_{X/Y}= B^{t,(0)}_{X/Y}$ (see \cite[Sec. 1.4,
p. 21]{bruns-herzog}, \cite[Ch. 20]{eisenbud:commutative} and
\cite{kallstrom:branchpurity}). If $X$ is not connected, then the
schemes $B^{(i)}_{X/Y}$ and $B^{t, (i)}_{X/Y}$ are defined on each
connected component of $X$, and therefore we have well-defined
subschemes of $X$.

We will study base change diagrams
\begin{displaymath}\tag{$BC $}
    \xymatrix{
 X_1 \ar[r] ^j\ar[d]^{\pi_1}  
& X\ar[d]^\pi\\
      Y_1 \ar[r]^i & Y,
    }
  \end{displaymath}
  where $X_1=X\times_{Y}Y_1$. If $Y_1 \to Y$ is flat or
  $B_{X/Y}= \emptyset$, then the canonical homomorphism  
\begin{displaymath}
\psi: j^*(T_{X/Y})\to T_{X_1/Y_1}.
\end{displaymath}
need not be an isomorphism (but see
the proof of \Proposition{\ref{spec-der-lemma}} below), contrary
to the good behaviour of $\Omega_{X/Y}$. On the other hand,
\Theorem{\ref{Z-L}} gives a rather general situation when $B_{X/Y}^t$
restricts nicely to the fibres of $X/Y$.

We start with the following important fact (see e.g. \cite[Proposition
16.4.5]{EGA4:4} for (1) below).

\begin{prop} \label{base-change-differentials} Consider the diagram
  (BC).
  \begin{enumerate}\item 
    The canonical morphism
    \begin{displaymath}
      j^*(\Omega_{X/Y})\to \Omega_{X_1/Y_1}
    \end{displaymath}
    is an isomorphism.
  \item Let $X'_1$ be a connected component of $X_1$.  Then
    $B^{(i)}_{X'_1/Y_1}= B^{(i+r)}_{X/Y}\times_X X'_1$, where $r=
    d_{X'_1/Y_1}- d_{X/Y}$.
\item   
If $\psi $ is an isomorphism, then
  \begin{displaymath}
    B^{t,r}_{X'_1/Y_1} =     B_{X/Y}^{t, r} \times_X X'_1,
  \end{displaymath}
  where $X'$ and $r$ are as in (2). In particular, if $\psi$ is an
  isomorphism, then $B^t_{X'_1/Y_1}= B_{X/Y}^t\times_XX'_1$.
\end{enumerate}
\end{prop}

\begin{proof} For the proof we can assume that all schemes are affine,
  so let $A\to B$ and $A\to A_1$ be homomorphisms of commutative
  rings

  (1): Put $B_1=A_1\otimes_A B$.  Let $d_{B/A}: B \to \Omega_{B/A}$ be
  a universal derivation and define the $A_1$-linear derivation $d=
  \id \otimes d_{B/A}: B_1 \to A_1\otimes_A \Omega_{B/A}$, which can
  be factorised over a universal derivation $d_{B_1/A_1}:B_1 \to
  \Omega_{B_1/A_1}$ by a $B_1$-homomorphism $\tilde d:
  \Omega_{B_1/A_1} \to A_1 \otimes_A \Omega_{B/A}$. There exists a
  natural $B_1$-homomorphism $p: A_1 \otimes_A \Omega_{B/A} \to
  \Omega_{B_1/A_1} $, which is the inverse of $\tilde d$.

  (2): Let now $B_1$ the affine ring of a connected component $X_1'$
  of $\Spec (A_1 \otimes_A B)$. By (1),
  $\Omega_{B_1/A_1} \cong B_1 \otimes_B \Omega_{B/A}$. Let
  $F_{i+d_{X_1/Y_1}}(\Omega_{B/A})$ denote the Fitting ideal defining
  $B^{(i+r)}_{B/A}$ and recall that for any Fitting ideal and
  $B$-module $M$ of finite type, $B_1 F(M)= F(B_1\otimes_B M)$. Then
  \begin{eqnarray*}
    B^{(i+r)}_{X/Y}\times _X X'_1 &=& V(F_{i+d_{X'_1/Y_1}}(\Omega_{B/A})) \times_{\Spec B}
    \Spec B_1 \\  &=& \Spec (\frac {B}{F_{i+d_{X'_1/Y_1}}(\Omega_{B/A})} \otimes_ B B_1 )  
    = \Spec
    \frac {B_1}{B_1 F_{i+d_{X'_1/Y_1}}(\Omega_{B/A})}\\ &=& \Spec \frac
    {B_1}{F_{i+d_{X'_1/Y_1}}(B_1\otimes_B \Omega_{B/A})} = \Spec \frac
    {B_1}{F_{i+d_{X'_1/Y_1}}(\Omega_{B_1/A_1})}\\ &=& B^{(i)}_{X'_1/Y_1}.
\end{eqnarray*}

(3) is proven in the same way as (2).
\end{proof}

\begin{proposition}\label{spec-der-lemma}
  Let $\pi: X/S\to Y/S$ be a finitely presented morphism of schemes,
  such that $B_{X/S}= \emptyset$. Put $\Ac = \supp \Ker (\psi)$.
  \begin{enumerate}
  \item Assume that
    \begin{enumerate}
    \item $X$ is Cohen-Macaulay and connected and $\pi$ is flat of
      relative dimension $d_{X/Y}>0$.
    \item $Y_1\to Y$ is a locally complete intersection morphism
    \end{enumerate}
    Then $\codim^+_{X_1} \Ac =0$. If
    $\psi$ is injective, then $\psi$ is an isomorphism.
  \item Make the assumption (a)  in (1). Let $X_y$ be the fibre over
    a point $y$ such that $\Oc_{Y,y}$ is regular, and $\pi$ is smooth
    at the generic points of $j(X_y)\subset X$. Then $\psi_y $ is an
    isomorphism.
 \end{enumerate}
\end{proposition}
 \begin{lemma}\label{djuplemma} 
      \begin{enumerate}\item
        Let $B$ be a ring, $I$ an ideal, and $N$ and $M$ be
        $B$-modules (not necessarily of finite type).  Any $N$-
        sequence in $I$ of length $2$ is also a $Hom_B(M,N)$ -
        sequence.

        In particular, any $B$-sequence of length $\leq 2$ is also a $
        T_{B/A}$-sequence, for any subring $A\subset B$.
      \item Let $(A,\mf_A) \to (B,\mf_B)$ be a flat homomorphism of
        local rings, and let $I\subset A$ be an ideal which is
        generated by an $A$-sequence. Let $N$ be a $B$-module of
        finite type which is flat over $A$.  If $\depth_{\mf_B} N/ I N
        \geq 1$ and $\depth_{\mf_B} N \geq 2$, then
        $\depth_{\mf_B} I Hom_B(M,N)\geq 2$.

      \item Let $(A,\mf_A) \to (B,\mf_B)$ be a flat homomorphism of
        local rings, where $A$ is regular. Let $N$ be a $B$-module of
        finite type which is flat over $A$.  If $\depth_{\mf_B}
        N/\mf_A N \geq 1$ and $\depth_{\mf_B} N \geq 2$, then
        \begin{displaymath}
        \depth_{\mf_B}\mf_AHom_B(M,N)\geq 2.
      \end{displaymath}
\end{enumerate}
\end{lemma}

\begin{pf} (1): (this is of course well-known) Let $x_1,x_2$ be an
  $N$-sequence in $I$. Clearly, $x_1$ is $Hom_B(M,N)$-regular. Assume
  $x_2\phi_2(m)= x_1\phi_1(m)$, $\phi_i \in Hom_B(M,N)$, $m\in M$.
  Since $x_1,x_2$ is an $N$-sequence, $\phi_2(m)= x_1 n'$, $n'\in N$.
  Since $x_1$ is a regular element this gives a well-defined
  homomorphism $\phi'\in Hom_B(M,N)$, $\phi'(m)= n'$, and
  $\phi_2=x_1 \phi' $, hence $x_1,x_2$ is a $Hom_B(M,N)$-sequence.

  (2): Let $x_1$ be $N/ I N$-regular, hence by flatness it is
  $N$-regular, and we can find $x_2$ such that $x_1,x_2$ forms an
  $N$-sequence. By assumption, $I= (y_1, \dots , y_r)$ where
  $y_1, \dots , y_r$ is an $A$-sequence, and since $N$ is flat it is
  also an $N$-sequence. Then $y_1, \dots , y_r, x_1$ is an
  $N$-sequence (see proof of \cite{matsumura}*{Th. 23.3}). Assume
  $x_2\phi_2 = x_1 \phi_1$, where $\phi_1, \phi_2 \in I Hom_B(M,N)$.
  As $x_1,x_2$ is $N$-regular, $\phi_2 (m)\in x_1N$, and since $x_1$
  is $N$-regular $\phi_2 = x_1\phi_2' $, where
  $\phi'_2 \in Hom_B(M,N)$. Therefore
  $\phi_2 \in x_1 Hom_B(M,N)\cap I Hom_B(M,N)$. Assume that
  $\sum y_i f_i = x_1 f$ where $f,f_i\in Hom_B(M,N) $. Since
  $y_1, \dots , y_r, x_1$ is an $N$-sequence we have
  $x_1 N\cap I N = x_1 IN$, hence $f_i(m)\in x_1 N$, and since $x_1$
  is $N$-regular, $f_i = x_1 f'_i$ where $f'_i \in Hom_B(M,N)$. This
  implies
  $\phi_2 \in x_1 Hom_B(M,N)\cap I Hom_B(M,N)= x_1 IHom_B(M,N)$, and
  thus $x_1, x_2$ is a $ I Hom_B(M,N)$-sequence.

  (3): This follows from (2) since $\mf_A$ is generated by an
  $A$-sequence.
\end{pf}

\begin{pfof}{\Proposition{\ref{spec-der-lemma}}}
$  (1)\Rightarrow (2)$: Since $\Oc_{Y,y}$ is regular, $\mf_{Y,y}$ is
  generated by an $\Oc_{Y,y}$-sequence, so that (1) is applicable.
  Since $d\pi$ is smooth at the generic points of $X_y$ it follows
  that $\Ac = \emptyset$.

  (1): By assumption $Y_1 \to Y$ can be locally factorised
  $i= i_2 \circ i_1$, where $i_1: Y_1 \to Y_2$ is a regular immersion
  (see \cite{EGA4:4}*{Def. 16.9}) and $i_2 : Y_2 \to Y$ is flat; the
  assertion being local on $Y_1$ we can in fact assume that we have a
  globally defined factorisation.

  Let $X_2 \to Y_2$ be the base change of $\pi$ over $i_1$,
  $j_1: X_1 \to X_2 $ be the second projection, and $X_1 \to Y_1$ the
  base change of $X_2/Y_2$ over $i_1$, and $j_2: X_2 \to X_1$ be the
  second projection. Since $\Omega_{X/Y}$ is of finite presentation,
  it follows from \cite{matsumura}*{Th. 7.11} that
  $j_1^*(T_{X/Y}) = T_{X_2/Y_2}$, and thus
  $j^*(T_{X/Y})= j_2^* (j_1^*(T_{X/Y})) = j^*_2(T_{X_2/Y_2})$. It is
  straightforward to check also that $d_{X/Y}> 0$ will remain true
  after base change. We can therefore reduce to the case when
  $i: Y_1 \to Y$ is a regular immersion. Let $I$ be the ideal of
  $j(X_1)$ in $X$.

  Consider the exact sequence
\begin{displaymath}
  0 \to T_{X/Y}\to T_{X/S} \xrightarrow{d\pi}  T_{Y/S \to X/S},
\end{displaymath}
where $T_{Y/S \to X/S} $ is the sheaf of $S$-derivations
$\pi^{-1}(\Oc_Y)\to \Oc_X$, and the tangent map
\begin{displaymath}
  T_{X_1/Y_1}\xrightarrow{dj} T_{X/S \to X_1/S } \xleftarrow{\nu}  j^*(T_{X/S}),
\end{displaymath}
where the map $\nu$ is an isomorphism since $\Omega_{X/S}$ is locally
free ($B_{X/S}= \emptyset$). We have therefore the diagram
    \begin{displaymath}
\label{galois-diagram}
\begin{tikzcd}
    T_{X/Y} \arrow[r] \arrow[d, "s"] & T_{X/S}
    \arrow[d]\arrow[r,"d\pi"]& T_{Y/S\to X/S}\arrow[d] \\
    j_*j^*(T_{X/Y}) \arrow[r,"\alpha"]\arrow[d, "\bar \psi"]&
    j_*j^*(T_{X/S})\arrow[r,"\overline{d\pi}"]\arrow[d,"\bar \nu"]& j_*j^*(T_{Y/S\to X/S}) \\
    j_*(T_{X_1/Y_1})\arrow[r,"dj"] & j_*(T_{X/S\to X_1/S}), &
\end{tikzcd}
  \end{displaymath}
  where $\bar \nu$ is an isomorphism and
  $\alpha = \bar \nu^{-1}\circ dj\circ \bar \psi$. Put
  $\psi_1= \bar \psi\circ s$ and $\Ac^1 = \Ker (\psi_1)$. Assume on the
  contrary that there exists a maximal point $x$ of $ \Ac \subset X_1$
  of height $\geq 1$. Let $x_1$ be a point in $X$ that specializes to
  $j(x)$, $x_1 \neq j(x)$, and $x_2$ be a maximal point in $X_1$ such
  that we have specializations
  \begin{displaymath}
   x_1 \rightsquigarrow j(x_2) \rightsquigarrow j(x).
 \end{displaymath}
 Then $x_2 \rightsquigarrow x$ and $x_2 \neq x$; hence
 $x_2\not \in \Ac$ and since $x_1 \rightsquigarrow j(x_2) $,
 $\Ac^1_{x_1} = I_{x_1} T_{X/Y,x_1}$. Letting
 $\phi: \Spec \Oc_{X,j(x)} \setminus \{\mf_{X,j(x)}\}\to \Spec
 \Oc_{X,j(x)}$ be the canonical open inclusion we thus have
 $\Ker (\psi_1)_{x_1} = I_{x_1}T_{X/Y,x_1}$ when $x_1$ is a point in
 $\Spec \Oc_{X,j(x)} \setminus \{\mf_{X,j(x)}\}$; hence
 $\phi^*(\Ker (\psi_1)) = \phi^*(I_{\pi(x)}T_{X/Y,x})$. Since $X$ is
 Cohen-Macaulay and $X/Y$ is flat of relative dimension $d_{X/Y}>0$,
 it follows that
 \begin{displaymath}
   \depth_{j(x)} \frac{\Oc_{X,j(x)}}{ I_{\pi(x),j(x)}} \geq 1.
\end{displaymath}
Therefore, by \Lemma{\ref{djuplemma}, (3),}
$\depth_{j(x)} I_{ \pi(x)} T_{X/Y,j(x)}$ $ \geq 2$, so that the local
cohomology
\begin{displaymath}
H^1_{\mf_{X,j(x)}}(I_{ \pi(x)} T_{X/Y,j(x)})=0.
\end{displaymath}
Since $T_{X/Y}$ is
torsion free, hence $\Ac^1_{j(x)}$ and $I_{ \pi(x)} T_{X/Y,j(x)}$  are torsion free, we get
\begin{displaymath}
  \Ac^1_{j(x)} \subset \phi_*\phi^*(\Ac^1_{j(x)}) = \phi_*\phi^* (I_{\pi(x)} T_{X/Y,j(x)}) = I_{\pi(x)} T_{X/Y,j(x)}
  \subset    \Ac^1_{j(x)}, 
\end{displaymath}
and therefore $\Ac^1_{j(x)} = I_{\pi(x)} T_{X/Y,j(x)}$. This
contradicts the assumption that $x\in \Ac$. The map $dj$ is injective,
$\nu$ is an isomorphism, and
$\Imo (\bar \nu^{-1}\circ dj)\subset \Imo (\alpha)$. Therefore, if
$\psi$ is injective, it is also an isomorphism.
\end{pfof}

\begin{lemma} \label{heightlemma} Assume that $x \rightsquigarrow x_1$
  and put $y= \pi(x)$ and $y_1= \pi(x_1)$.
  \begin{enumerate}
  \item If $x_1 $ is maximal in the fibre $X_{\pi(x_1)}$ with this
    property (so that if
    $x \rightsquigarrow x'_1 \rightsquigarrow x_1 $ and
    $\pi(x'_1)= y_1$, then $x'_1 = x_1$). Then
    \begin{displaymath}
      \hto_X(x_1/x) \leq  \hto_Y(y_1/y).
    \end{displaymath}
  \item Assume that $X$ is catenary and equidimensional at $x_1$ (see
    \Theorem{\ref{hartshorne}}) and that $\pi$ is flat at $x_1$. Then
  \begin{enumerate}
  \item $ \hto_{X_{y_1}}(x_1) \leq \hto_{X_{y}}(x).$
  \item $ \hto_{X_{y_1}}(x_1)\leq \hto_X(x).$
  \end{enumerate}
\end{enumerate}
\end{lemma}
\begin{proof}
  (1): Let $V $ and $W$ be closure of $x$ and $y$ in $X$ and $Y$,
  respectively, and put $A= \Oc_{W, y_1 }$ and $B= \Oc_{V, x_1}$, so
  that we have a homomorphism $A \to B$ of local noetherian rings, and
  $\hto_X(x_1/x) = \dim B$ and $\hto_Y(y_1/y) = \dim A$. Since
  $x_1 \in \Max (V \cap X_{y_1})$, we have $\dim B/ \mf_A B =0$, so
  that by the dimension inequality \cite{matsumura}*{Th 15.1}
  \begin{displaymath}
    \dim B \leq  \dim A + \dim   B/ \mf_A B = \dim A,
  \end{displaymath}
  which implies the assertion.

  (2): (a) Since $\pi$ is flat at $x_1$, by the dimension equality [loc
  cit] the assertion is equivalent to
  \begin{displaymath}
    \hto_X(x_1) - \hto_X(x) \leq \hto_Y(y_1) - \hto_Y(y).
  \end{displaymath}
  Since $X$ is catenary and equidimensional at $x_1$,
  $\hto_X(x_1) - \hto_X(x) = \hto(x_1/x)$, and since $Y$ is catenary
  and equidimensional at $y_1$ \cite{matsumura}*{Th. 31.5 } so that
  $ \hto_Y(y_1) - \hto_Y(y) = \hto_Y(y_1/y)$, the assertion follows
  from (1).

  (b) By flatness, the dimension equality implies
  \begin{displaymath}
    \hto_X(x) = \hto_{X_y}(x) + \hto_Y(y) \geq  \hto_{X_y}(x) \geq
    \hto_{X_{y_1}}(x_1),  
  \end{displaymath}
  by (a).
\end{proof}


It is proven in \cite{SGA2}* {Exp. II} that if a catenary scheme $X$
satisfies $\depth (x) \geq 2$ when $\hto_X(x) \geq 2$, then $X$ is
equidimensional. We give a short account of the main content in [loc
cit, Exp III, Cor. 3.7, Cor 3.9], based on Hartshorne's connectedness
theorem (no originality is claimed).

\begin{theorem}\label{hartshorne}
  \begin{enumerate}
  \item Let $X$ be locally noetherian and connected scheme such that
    the depth $\depth (x) \geq 2$ when $x$ is a point of height
    $\hto_X(x) \geq d$.
    \begin{enumerate}
      \item Given $x, x' \in \Max(X)$ there exists a sequence
      $x'= x_0 , x_1, \dots, x_n = x \in \Max (X)$ such that
      $x_i^- \cap x_{i+1}^- \neq \emptyset$.
      \item  If
    $ \xi \in \Max (x^-_i \cap x^-_{i+1})$, then
    $\hto_X(\xi)\leq d-1$. ($X$ is connected in codimension $d-1$)
    \end{enumerate}
  \item Let $X$ be the spectrum of a local noetherian catenary ring,
    and assume that $\depth (x) \geq 2$ when $\hto_X(x) \geq 2$. Then
    $X$ is equidimensional, i.e. if $x_1, x_2 \in \Max(X)$, then
    $\hto_X(x_c/x_1)= \hto _X(x_c/x_2) $, where $x_c$ is the closed
    point in $X$.
  \end{enumerate}
\end{theorem}
\begin{proof}
  (1): The existence of the connecting chain of maximal points as
  stated, for any two maximal points $x'$ and $x$, follows since $X$
  is connected and locally noetherian. To prove the other assertion we
  can restrict to an open neighbourhood of point in $\xi$
  $ \Max (x^- _i\cap x^-_{i+1})$, so that it is the sole member of the
  set $ \Max (x^- _i\cap x^-_{i+1})$. Then
  $(x^-_i \cup x^-_{i+1}) \setminus \xi^-$ is a disconnected set.
  According to Hartshorne's theorem [loc cit, Théorème 3.6] it follows
  that $\hto_X(\xi) \leq d-1$.

  (2): We use the same notation as in (1). Taking a sequence of
  maximal points in $X$ as in (1) it suffices to see that
  $\hto_X(x_c/x_i ) = \hto_X(x_c/x_{i+1})$. By (1),
  $\hto_X(\xi) \leq 1$, and since $\xi \not\in \Max(X)$, we have
  $\hto_X(\xi)= \hto_X(\xi/x_i)= \hto_X(\xi/x_{i+1})=1$. Since $X$ is
  catenary,
  $\hto_X(x_c/x_i) = \hto_X(x_c/\xi) + \hto_X(\xi/x_i) =
  \hto_X(x_c/\xi) + \hto_X(\xi) = \hto_X(x_c/x_i)$.
\end{proof}

\section{Proof of  \Theorem{\ref{Z-L}}}

\begin{proposition}\label{lipman-normal}
  \begin{enumerate}
  \item \label{lichtenbaum-schlessinger}(Lichtenbaum-Schlessinger
    \cite [Prop. 5.2]{lipman:freeder}) Let $X/Y$ be a d.c.i. morphism
    which is locally of finite type and $X$ is Cohen-Macaulay. If
    $T_{X/Y}$ is locally free, then $\codim_X^+ B_{X/Y} \leq 2$.
    \item (Lipman \cite[Th. 1]{lipman:freeder}) Let $X/k$ be a
  scheme locally of finite type over a field of characteristic $0$
  such that $B^t_{X/k}=\emptyset$. Then $X$ is normal, and in
  particular $\codim^-_X B_{X/k} \geq 2$.
  \end{enumerate}
\end{proposition}
We include a proof, following [loc cit], to clarify the
situation in our notation.
\begin{proof}
  (1): Since $\pdo \Omega_{X/Y,x}\leq 1$ for each point $x$ in $X$, it
  follows that $B_{X/Y} = \supp Ext^1_{\Oc_X}(\Omega_{X/Y}, \Oc_X)$.
  Let $0 \to F_1 \to F_0 \to \Omega_{X/Y,x}\to 0 $ be a free extension
  of the $\Oc_{X,x}$-module $\Omega_{X/Y,x}$, which after dualisation
  gives the exact sequence
    \begin{displaymath}
      0 \to T_{X/Y,x}\to F_0^* \to F_1^* \to Ext^1_{\Oc_{X,x}}
      (\Omega_{X/Y,x}, \Oc_{X,x}) \to 0. 
    \end{displaymath}
    Since $T_{X/Y,x}$ is free, it follows that $\pdo Ext^1_{\Oc_{X,x}}
    (\Omega_{X/Y,x}, \Oc_{X,x}) \leq 2$. If $x\in \Max(B_{X/ Y})$, so
    $\depth (Ext^1_{\Oc_{X,x}} (\Omega_{X/Y,x}, \Oc_{X,x}))=0$, by the
    Auslander-Buschsbaum formula, $\depth (x) = \pdo Ext^1_{\Oc_{X,x}}
    (\Omega_{X/Y,x}, \Oc_{X,x}) \leq 2$, hence, since $\Oc_{X,x}$ is
    Cohen-Macaulay, $\hto_X(x)\leq 2$.

    (2): One first proves that the module $ C_X:= \Coker
  (\Omega_{X/k}\xrightarrow{g} T^*_{X/k})$, where $g$ is the biduality
  morphism, satisfies $\depth C_{X,x}\geq 2$ when $x\in\supp C_X $.
  Consider the exact sequence $0\to \overline {\Omega}_{X/k}\to
  T^*_{X/k}\to C_X \to 0$, where $\overline {\Omega}_{X/k}=
  \Imo(\Omega_{X/k}\to T^*_{X/k})$.  Noting that
  $T_{X/k}=\Omega_{X/k}^* = \overline \Omega_{X/k}^*$, dualisation results in the
  exact sequence
    \begin{displaymath}
      0\to C^*_X\to    T_{X/k}^{**} \to T_{X/k} \to Ext^1_{\Oc_X}(C_X,\Oc_X) \to 0 
  \end{displaymath}
  since $T_{X/k}^*$ is locally free. As $T_{X/k}$ is reflexive we get
  \begin{displaymath}
  Ext^0_{\Oc_X}(C_X,\Oc_X)= C_X^*=0 \quad \text {and} \quad Ext^1_{\Oc_X}(C_X, \Oc_X)=0,
\end{displaymath}
implying the assertion.
 
We always have $\supp C_X \subset B_{X/k}$. If $x\not \in \supp C_X$,
so the map $\Omega_{X/k,x}\to T_{X/k,x}^*$ is surjective, by a result
of Nagata \cite{nagata:remark-zariski} there exist
$\partial_i\in T_{X/k,x}$ and $x_j\in \mf_{X,x}$ such that
$\partial_i(x_j)$ forms an invertible $d\times d$ matrix, where
$d= \hto (x)$. Since $\Char k=0$ it follows from the
Zariski-Lipman-Nagata criterion that $\Oc_{X,x}$ is a regular ring,
hence again since $\Char k =0$, $x\notin B_{X/k}$. Therefore
$\supp C_X = B_{X/k}$. By the first assertion in the Proposition, it
follows that $\depth (x)\geq 2$ when $x\in B_{X/k}$. Since regularity
implies normality, so the locus of points where $X$ fails to be normal
is contained in $B_{X/k}$, it follows that $X$ is normal (either look
at Lipman's nice argument in [loc cit, Prop 2.1] or think of Serre's
normality criterion).
\end{proof}
\begin{pfof}{\Theorem{\ref{Z-L}}} 
(1):  Assume that $x \not\in B^t_{X/Y}$ so that $T_{X/Y,x}$ is free and
  select $\xi \in \Max (B_{X_y/k_{Y,y}})$, where $y=\pi(x)$, such that
  $ \xi \rightsquigarrow x$, so that $y=\pi(\xi)= \pi(x)$. There
  exists a point $x_m\in \Max(B_{X/Y})$ that specialises to $\xi$; put
  $y_m= \pi(x_m)$. Since $X$ is Cohen-Macaulay, or more generally
  $\depth (x')\geq 2$ when $\hto_X(x')\geq 2$ for any point $x'$ in
  $ X$, by \Theorem{\ref{hartshorne}} $X$ is equidimensional, so by
  \Lemma{\ref{heightlemma}} (2)
\begin{displaymath}
  \hto_{X_{y}}(\xi) \leq \hto_{X_{y_m}}(x_m). 
\end{displaymath}
Since $X$ is Cohen-Macaulay and $X/Y$ is a d.c.i, $\hto_X(x_m)\leq 2$
\Prop{\ref{lipman-normal}(1)}, and since
$\codim_Y^- D_\pi \geq 1$, $\hto_Y(y_m)\geq 1$. Therefore, by the
dimension equality for flat morphisms
  \begin{displaymath}
    \hto_{X_{y_m}}(x_m) = \hto_X(x_m) - \hto_Y(y_m) \leq 2-1 = 1.
  \end{displaymath}
  Therefore $\hto_{X_{y}}(\xi)\leq 1$.

  (2): Since $X_y/k_{k_{Y,y}}$ is generically smooth and $\Oc_{Y,y}$
  is regular, \Proposition{\ref{spec-der-lemma}}(2) implies that
  $\psi: j^*(T_{X/Y})\cong T_{X_y}$, hence by
  \Proposition{\ref{base-change-differentials}}(3),
  $ B^t_{X/Y} = j(X_y)\cap B_{X/Y}$.
\end{pfof}
\begin{pfof}{\Corollary{\ref{Z-L-cor}}} Put $A= k[y_1, \dots , y_r] $ and
  $B =k[X_1, \dots , X_n]$. If $\{f_1, \dots , f_r \}$ is a regular
  sequence in $ k[X_1, \dots , X_n]$ defining $V$, then $V=X_y$ is a
  fibre of the flat morphism $\pi: X= \Spec B \to \Spec A$,
  $y_i \mapsto f_i$ (see \cite[Exercise 22.2]{matsumura}), where $y$
  denotes the maximal ideal $(y_1, \dots , y_r)\subset A$. Since
  $T_{X_y/k_{Y,y}}$ is locally free and by assumption
  $\codim^-_Y D_{X/Y}\geq 1$, \Theorem{\ref{Z-L}}(2) implies that
  $x\notin B^t_{X/Y} $, so that by (1)
  $\codim^+_V (B_{V/k}) = \codim^+_{X_y}(B_{X_y/k_{Y,y}})\leq 1$. If
  moreover $\Char k =0$, \Proposition{\ref{lipman-normal}}(2) implies
  that $B_{V/k}= \emptyset$.
\end{pfof}

\subsection*{Hypersurfaces}
The proof of \Corollary{\ref{Z-L-cor}} is not conducted by reducing to
the case of hypersurfaces, and unlike Scheja and Storch's proof
\cite{scheja-storch:eigenschaften}*{Satz 9.3} for hypersurfaces $X/k$
over fields $k$ of characteristic $0$ the Eagon-Northcott bound on
heights of determinantal ideals is not utilized. Since their proof in
the hypersurface case is quite different it can be instructive to
compare proofs.

The proof below a little more geometric but still quite similar to
that in [loc cit]. In fact, Scheja and Storch even allow noetherian
$k$-algebras that can be provided with a universal finite derivation;
it is likely that the proof of \Corollary{\ref{Z-L-cor}} also can be
generalizes to such algebras.

\begin{proposition}\label{scheja-storch}
  Assume that $j(X)$ is a hypersurface in $Z$. Then
  $B^t_{X/k}= B_{X/k}$. In particular, if $T_{X/k }$ is locally free,
  then $X/k$ is smooth.
\end{proposition}
\begin{proof} We assert that if $x\in \Max(\Bc_{X/k})$, then
  $\hto_X(x)\leq 1$. Assume on the contrary that $\hto_X(x)\geq 2$.
  Put $\bar B= \Oc_{X,x}$ and $B= \Oc_{Z,j(x)}$, so that by assumption
  $\bar B= B/I$, where $I= (f)$ for some element $f$ in the smooth
  local $k$-algebra $ B$. We have the short exact sequence
  $0\to I/I^2 \to \bar B\otimes_B \Omega_{B/k}\to \Omega_{\bar B/k}\to
  0$, where $I/I^2\cong \bar B$. Dualising we get the exact sequence
\begin{displaymath}
  0\to T_{\bar B/k} \to \bar B\otimes_B(T_{B/k})\to \bar \Cc \to 0, \tag{E} 
\end{displaymath}
Since $T_{\bar B/k}= \bar B \otimes T_{B/k}(I)$, where $T_{B/k}(I)$ is
the submodule of derivations such that $\partial (f)\in (f)$, it
follows that
\begin{displaymath}
  \bar \Cc \cong \frac {T_{B/k}\cdot f}{I} \subset \bar B, 
\end{displaymath}
so that $\bar \Cc$ can be identified with an ideal of $\bar B$. Notice
that
\begin{displaymath}\tag{*}
\bar \Cc_{x'} = \bar B_{x'} 
\end{displaymath}
when $x' $ is a point in $\Spec \bar B$ of height $\leq 1$. It follows
that $\bar \Cc$ contains a regular sequence of length $\geq 2$, since
$\bar B$ is normal. This can also be concluded from the existence of
\thetag{E}, since by the Hilbert-Burch theorem it implies that
$\bar \Cc = a F_1(\bar \Cc)$ (see \cite[Th.
20.15]{eisenbud:commutative}) for some non-zero divisor $a$, and
$F_1(\bar \Cc)$ has depth exactly $2$. By \thetag{*} and Krull's
theorem it follows that $a$ is a unit, so that ni fact
$\bar \Cc = F_1(\bar \Cc)$. One can extend \thetag{E} to an exact
sequence of $B$-modules
\begin{displaymath}\tag{E'}
0 \to   B^d \to B^{d+1} \to \Cc \to 0,
\end{displaymath}
where $\Cc = T_B\cdot f + (f)$. Now $f\in (T_B\cdot f)_{x'}$ when $x'$
is a point of height $\leq 1$, hence
\begin{displaymath}
  f \in \bigcap_{x', \hto(x')\leq 1}  (T_B\cdot f)_{x'} =
  \sqrt{T_B\cdot f},
\end{displaymath}
the integral closure if $T_{B/k}\cdot f$; therefore 
\begin{displaymath}
  \Max_{\Spec \bar B}(V(\bar \Cc))= \Max_{\Spec B } (V(\Cc)).
\end{displaymath}
By the Eagon-Nortcott bound \cite{eagon-northcott:height}, applied to
\thetag{E'}, any point in the above set is of height $\leq 2$ in
$\Spec B$, hence it is of height $\leq 1$ in $\Spec \bar B$. This
contradicts \thetag{*}. Therefore if $\Max (B_{X/k})\neq \emptyset$ it
follows that $\hto_X(x)\leq 1$, but since $X$ is normal, it follows
that $B_{X/k}= \emptyset$.

\end{proof}
\section{ Differential criterion of smoothness}\label{sect:2}
The relation between the ramification locus and the locus of
non-smooth points of a morphism is of course much discussed in the
literature, but there still seems to remain room for clarification. In
\cite[\S 2]{spivakovsky:popescu} one can find a nice summary of
characterisations of smoothness in terms of the vanishing of
Andr{\'e}-Quillen homology and the Jacobian condition
$B_\pi= \emptyset$. We are however more interested in the ``Jacobian
characterisation'' \cite[Proposition 17.15.15]{EGA4:4}. A proof of the
relevant statement (\Theorem{\ref{smooth-flat}}) is included since I
find the argument in EGA difficult to disentangle and was unable to
find any other satisfactory treatment of this important result. The
proof below relies on the assertions that smoothness implies flatness,
and conversely, if a morphism is flat at a point and smooth along the
fibre at the point, then the morphism is smooth at that point
\cite[Théorème 19.7.1]{EGA4:1}.

\begin{thm}\label{smooth-flat}
  Let $\pi: X\to Y$ be a morphism of schemes which is locally finitely
  presented. Let $x$ be a point in $X$ and put $y=\pi(x)$. The
  following are equivalent:
  \begin{enumerate}[label=(\theenumi),ref=(\theenumi)]
  \item\label{fl-smooth1} $\pi$ is smooth at the point $x$.
  \item\label{fl-smooth2} $\pi$ is flat at $x$, $x\not\in B_\pi$, and
    $X_y/k_{Y, y}$ is smooth at all points $\xi \in \Max
    (X_y) $ that specialise to $x$.
  \item\label{fl-smooth3} $\pi$ is flat at $x$, $x\not\in B_\pi$, and
    $\rko \Omega_{X/Y,x}= \dim_x X_y$.
  \item\label{fl-smooth4} $\pi$ is flat at $x$, $x\not\in B_\pi$, and
    $X_y/k_{Y,y}$ is smooth at all points $\eta \in \Max (X_y) $ that
    specialise to $x$.
  \end{enumerate}
\end{thm}
\begin{remark}
  The proof in \cite{EGA4:4} that the rank of $\Omega_{X/Y,x}$ is as
  asserted in $(1)\Rightarrow (3)$ seems incomplete. It relies on
  \cite[Proposition 17.10.2]{EGA4:4}, and presupposes that $\pi$ be
  smooth not only at the point $x$, but also at all specialisations of
  $x$, in order to reduce to $k$-rational points and $Y= \Spec k$.
\end{remark}

\begin{proof}
  $ \ref{fl-smooth1} \Rightarrow \ref{fl-smooth2} \text{ and }
  \ref{fl-smooth3}$: Put $A=\Oc_{Y,y}$ and $B= \Oc_{X,x}$, so $A\to B$
  is formally smooth, and since moreover this homomorphism is finitely
  presesented, it is also flat \cite[Th 17.5.1]{EGA4:1}. Let $\xi$ be
  a point that specialises to $x$, and put $B'= \Oc_{X,\xi}$ and $A'=
  \Oc_{Y,\pi(\xi)}$. It follows directly from the definition of formal
  smoothness that the composition $A \to B \to B'$ and the base change
  $A' \to B'$ are also formally smooth (see also \cite[Theorem
  28.2]{matsumura}); hence $\pi$ is smooth at $\xi$. Since $A\to B$ is
  formally smooth it follows that $\Omega_{B/A}$ is projective
  \cite[Theorem 28.5]{matsumura}, and since $B$ is a local ring,
  $\Omega_{B/A}$ is free; hence $x\notin B_\pi$. We now determine the
  rank $r= \rank \Omega_{B/A}$. Let $k$ be an algebraic closure of
  $k_{Y,y}$, put $\hat B = k\otimes_{k_{Y,y}} k_{Y,y}\otimes_{A} B$,
  and let $k_{\hat B}$ be the residue field of $\hat B$.
  \Proposition{\ref{base-change-differentials}} implies that $
  \Omega_{\hat B/k}$ is free of rank $r$, hence by the second
  fundamental exact sequence in \cite[Theorem 25.2]{matsumura}, noting
  that $k_{\hat B}/k$ is formally smooth since $k$ is algebraically closed,
  \begin{displaymath}
    r =  \dim_{k_{\hat B}} \frac {\mf_{\hat B}}{\mf^2_{\hat B}} +
    \trdeg k_{\hat B}/k.
  \end{displaymath}
  Formal smoothness is preserved under base change, hence the map
  $k \to \hat B$ is formally smooth; hence $\hat B$ is a regular local
  ring \cite[Lemma 1]{matsumura}, so that
  $\dim_{\hat B} \frac {\mf_{\hat B}}{\mf^2_{\hat B}} = \dim \hat B =
  \dim B =\hto_{X_y}(x)$. Adding up we get
  \begin{displaymath}
    \dim_x X_y = \hto_{X_y}(x) + \hto_{X_y}(x_1/x)=r
\end{displaymath}
if we prove that $\hto_{X_y}(x_1/x) = \trdeg k_{\bar B}/k$ when $x_1$
is a closed point in $X_y$ such that $x \rightsquigarrow  x_1$. Let
$P_x = \mf_{X_y,x }^c$ be the contration of $\mf_{X_y,x }$ in
$\Oc_{X_y,x_1}$. We have (as detailed below)
  \begin{eqnarray*}
    \hto_{X_y}(x_1/x)& =& \dim \Oc_{X_y,   x_1}/ P_x= \trdeg
    K(\Oc_{X_y, x_1}/P_x)/k_{Y,y}\\ & =& \trdeg
    k_{X_y,x}/k_{Y,y} =  \trdeg k_{\hat B}/k.
  \end{eqnarray*}
  The first equality is by definition of Krull dimension, and the
  second, where $ K(\Oc_{X_y, x_1}/P_x)$ is the fraction field of
  $\Oc_{X_y, x_1}/ P_x$, follows since the integral domain $\Oc_{X_y,
    x_1}/ P_x$ is essentially of finite type over $k_{Y,y}$
  \cite{matsumura}*{Th 5.6}. The third equality is clear, and the last
  equality follows since extending $k_{X_y,x}$ and $k_{Y,y}$ by the
  algebraic closure of $k_{Y,y}$ will not change the transcendence
  degree.


  $ \ref{fl-smooth2} \Rightarrow \ref{fl-smooth1}$: Since $\pi$ is
  flat at $x$, it suffices by \cite[Th 19.7.1]{EGA4:1} to prove that
  $k_{Y,y}\to k_{Y,y}\otimes_{\Oc_{Y,y}}\Oc_{X_y, x}$ is formally
  smooth to conclude that $\Oc_{Y,y}\to \Oc_{X,x}$ is formally smooth,
  and by \cite[Theorem 30.3]{matsumura} (see also \cite[Proposition
  15.15.5]{EGA4:4}) this follows if $\Omega_{\Oc_{X_y,x}/k_{Y,y}}$ is
  free of rank $\dim_x X_y$.

  Since $x\notin B_\pi $, by
  \Proposition{\ref{base-change-differentials}} the
  $\Oc_{X_y,x}$-module $\Omega_{X_y/k_{Y,y},x} $ is free and its rank
  $r= \rko \Omega_{X_y/k_{Y,y},x} = \rko \Omega_{X/Y,x}$. Since $X_y
  \to \Spec k_{Y,y}$ is smooth at points $\xi \in \Max (X_y)$ that
  specialise to $x$, and $\Omega_{X_y/k_{Y,y}, \xi} =
  \Oc_{X_y,\xi}\otimes_{\Oc_{X,x}}\Omega_{X_y/k_{Y,y},x}$, it follows
  from $(1) \Rightarrow (2)$ (for the morphism $X_y \to \Spec
  k_{Y,y}$) that $r= \dim_\xi X_y $ (or see \cite[Prop.
  17.15.5]{EGA4:4}).  Since $X_y$ is locally of finite type over
  $k_{Y,y}$, the dimension theorem implies
  \begin{displaymath}
    \dim_x X_y = \trdeg k_{X_y,x}/k_{Y,y} + \hto_{X_y}(x)  =  \trdeg
    k_{X_y,\xi}/k_{Y,y} = \dim_\xi X_y =r.
  \end{displaymath}

  $ \ref{fl-smooth3} \Rightarrow \ref{fl-smooth4}$: Since
  $\Omega_{X/Y,x}$ is free of rank $\dim_x X_y $ it follows that
  $\Omega_{X_y/k_{Y,y},x}$ is free of rank $\dim_x X_y $
  \Prop{\ref{base-change-differentials}}.  In the same way as in the
  proof of $ \ref{fl-smooth2} \Rightarrow \ref{fl-smooth1}$ it follows
  that $k_{Y,y}\to \Oc_{X_y,x}$ is formally smooth; hence any
  localisation $k_{Y,y}\to \Oc_{X_y,\eta}$ is also formally smooth.

  $ \ref{fl-smooth4} \Rightarrow \ref{fl-smooth1}$: Since $k_{Y,y}\to
  \Oc_{X_y, \eta}$ is formally smooth it follows that
  $\Omega_{X_y/k_{Y,y},\eta}$ is free of rank $\dim_\eta X_y = \dim_x
  X_y$. Now $\Omega_{X/Y,x}$ is free and by
  \Proposition{\ref{base-change-differentials}} its rank is $\dim _x
  X_y$.  It follows as in the proof of $ \ref{fl-smooth2} \Rightarrow
  \ref{fl-smooth1}$ that the homomorphism $\Oc_{Y,y}\to \Oc_{X,x}$ is formally smooth.
\end{proof}
Often the condition \ref{fl-smooth2} or \ref{fl-smooth3} in
\Theorem{\ref{smooth-flat}} serve as a {\it definition} of smoothness
(see e.g. \cite{hartshorne}). Alternatively, $\pi$ is smooth if it is
flat and all its fibres are smooth \cite[Théorème 17.5.1]{EGA4:4}. In
either case, the condition that $\pi$ be flat can be a nuisance.  Put
$\Gamma_{X/Y/S}= \Ker (\pi^*(\Omega_{Y/S})\to \Omega_{X/S})$. Assuming
$X/S$ is smooth, \ref{smooth3} in the following theorem shows that the
non-smoothness locus of $\pi$ is exactly $\supp \Gamma_{X/Y/S}\cup
B_\pi$. Therefore, if $\Gamma_{X/Y/S}=0$ it follows that the Jacobian
criterion $B_\pi = \emptyset$ implies smoothness, i.e.  flatness is
automatic. Moreover, $\Gamma_{X/Y/S}=0$ when either $X/Y$ is
generically smooth and a locally complete intersection, or $X/S$ is a
locally complete intersection (see \cite[Proposition
2.11]{kallstrom:branchpurity} for a discussion of this assertion).



\begin{thm}\label{th-smooth} Assume that $X/S$ is smooth at the point $x$.
  The following are equivalent:
  \begin{enumerate}[label=(\theenumi),ref=(\theenumi)]
  \item\label{smooth1} $\pi: X \to Y $ is smooth at $x$.
  \item\label{smooth2} $\pi^*(\Omega_{Y/S})_x\xrightarrow{p}
    \Omega_{X/S,x}$ has a left inverse.
  \item\label{smooth3} $x\notin B_\pi $ and $\Gamma_{X/Y/S,x}=0$.
  \end{enumerate}
\end{thm}

\begin{proof}
  $ \ref{smooth1}\Leftrightarrow\ref{smooth2}$: See \cite[Th.
  17.11.1]{EGA4:4}.  $\ref{smooth3}\Rightarrow\ref{smooth2}$: This can
  be seen directly from the fundamental exact sequence of
  differentials \cite[Proposition 8.11]{hartshorne}.  $ \ref{smooth1}\Rightarrow\ref{smooth3}$:
  \ref{smooth1} implies by \Theorem{\ref{smooth-flat}} that
  $x\notin B_\pi$, and since \ref{smooth1} implies \ref{smooth2} we
  also get $\Gamma_{X/Y/S,x}=0$.
\end{proof}
\section{ Freeness and submersions}
This section contains some assorted results that can be regarded to be
related to \Theorem{\ref{Z-L}} but are not used for its proof.
\subsection*{Free modules}
We recall an easy extension of \cite[Chapter II, Lemma 8.9
]{hartshorne}:
 \begin{description}
\item[($*$)] a coherent $\Oc_X$-module $M$ is free at a point $x$ if $M_{\xi}$
 is free of rank equal to $ \dim_{k_{X,x}} k_{X,x}\otimes_{\Oc_{X,x}}M_x $ for
 each $\xi \in \Max (X)$ that specialises to $x$.
\end{description}

\begin{proposition}\label{freelemma} Let $M$ be a coherent $\Oc_X$-module and $x$ be a point
  in $X_1$, in the diagram (BC). Assume that $\depth (j(x))\geq 2$ and
  that $M_\xi$ is free for each point $\xi$ that specialises to $x$
  and $\xi \neq x$. Assume that there exists a point $x_1 $ in $X_1$
  that specialises to $x$ and $j(x_1)\neq j(x)$. The following are
  equivalent:
  \begin{enumerate}
  \item $j^*(M)_x$ is free. 
  \item $M_{j(x)}$ is free.
  \end{enumerate}
\end{proposition}
\begin{proof}
  Here $(2)\Rightarrow (1)$ is obvious, so assume (1).  By assumption
  there exists a neighbourhood $U $ of $j(x)$ such that, putting $U'=
  U \setminus \{j(x)\}$, the restriction $M_{U'}$ is locally
  free. Since $\depth (j(x))\geq 2$, Hartshorne's theorem \cite{SGA2}*
  {Th. 3.6} implies that $U'$ is connected, so the rank $r$ is the
  same at each point $\xi\neq x$ that specialises to $x$. Since
  $j(x_1)\neq j(x)$ specialises to $j(x)$, so $j(x_1)\in
  U'$. Therefore the rank of the free module $j^*(M)_x$ also equals
  $r$.  Since $k_{X,j(x)}\otimes_{\Oc_X,j(x)}M_{j(x)} = k_{X_1,
    x}\otimes_{\Oc_{X_1}, x} j^*(M)_x$, it follows from $(*)$ that
  $M_{j(x)}$ is free.
\end{proof}

\subsection*{Submersive non-smooth maps}
Let $\pi: X/S \to Y/S$ be a morphism of schemes that is locally finite
type, with ramificaiton locus $B_{X/Y}$ and discriminant locus
$D_{X/Y}= \pi(B_{X/Y})$. Let $T_{X/S\to Y/S}$ be the sheaf of
derivations $\pi^{-1}(\Oc_Y)\to \Oc_X$, so we have the exact sequence
\begin{displaymath}\tag{*}
  0 \to T_{X/Y}\to T_{X/S} \xrightarrow{d\pi} T_{X/S\to Y/S}
\end{displaymath}
where $d\pi$ is the tangent mapping.

The following result can be regarded as a weak converse of
\Proposition{\ref{lipman-normal}}(1). It shows also that all sections
of $T_{Y/S}$ can be locally lifted to sections of $T_{X/S}$ under
rather weak conditions on $\pi$.

\begin{theorem}\label{lemmasurjective}Assume:
  \begin{enumerate}
  \item $\pi$ is flat and $X$ is Cohen-Macaulay.
      \item $B_{X/S}= \emptyset$ (e.g. $X/S$ is smooth).
      \item $\codim^-_X B_{X/Y} \geq 3$.
  \item The base change $X_{D_{X/Y}}\to D_{X/Y}$ to the discriminant
    locus is generically smooth (i.e. generically reduced and
    separable)  and $D_{X/Y}$ is a local
        complete intersection in $Y$.
    \item $T_{X/Y}$ is $Y$-flat. This holds in particular when
  $\dim Y \leq 2$.
  \end{enumerate}
  Then $d\pi$ is surjective, and if $B^t_{Y/S}= \emptyset$, it follows
  that $B^t_{X/Y}= \emptyset$.
\end{theorem}
\begin{remark}
  If $X/Y$ is d.c.i. and $\pi$ is generically smooth, then
  $\codim^+ B_{X/Y} \leq d_{X/Y} +1$ (see e.g.
  \cite{kallstrom:branchpurity}).
\end{remark}
\begin{proof}
  We have the exact sequence
  \begin{displaymath}
  0 \to \Gamma_{X/Y/S} \to \pi^*(\Omega_{Y/S})\to \Omega_{X/S}\to
  \Omega_{X/Y}\to 0,
\end{displaymath}
where $\supp \Gamma_{X/Y/S}\subset B_{X/Y}$. Since also
$\supp Ext^1_{\Oc_X}(\Omega_{X/Y}, \Oc_X)\subset B_{X/Y}$, we get the
exact sequence
\begin{displaymath}\tag{**}
  0 \to i_*i^*T_{X/Y} \to i_*i^*(T_{X/S})\to  i_*i^*(T_{X/S \to Y/S}) \to R^1i_*i^*(T_{X/Y}),
\end{displaymath}
where $i: X_0 = X\setminus B_{X/Y} \to X$ is the open embedding.
We assert that: 
\begin{displaymath}\tag{\&}
  \depth_{B_{X/Y}} T_{X/Y} \geq 3, 
\end{displaymath}
and therefore
\begin{displaymath}
  R^1i_*i^*(T_{X/Y}) = R^2\Gamma_{B_{X/Y}} (T_{X/Y}) =0.
\end{displaymath}
Since $X$ and hence $T_{X/Y}$ and $T_{X/S\to Y/S}$ satisfy
$(S_2)$ \Lem{\ref{djuplemma}(1)}, one can  moreover erase $i_*i^*$ in the three first terms in
\thetag{**}. One can therefore add $\to 0$ to the right in \thetag{*}.

If $B_{Y/S}^t = \emptyset$, then \thetag{*} is locally split exact,
hence, since $B_{X/S}= \emptyset$ so that $B^t_{X/S}= \emptyset$, it
follows that $B_{X/Y}^t = \emptyset$.

We prove \thetag{\&}. By (5) $T_{X/Y}$ is $Y$-flat and $Y$ is
Cohen-Macaulay \cite{matsumura}*{\S 23 Corollary}, so that if  $x\in
B_{X/Y}$, by [loc. cit,
Th. 23.3] 
\begin{displaymath}\tag{\&\&}
  \depth_{x}(T_{X/Y,x})=  \hto_Y(y) +    \depth_{x} (j^*(T_{X/Y})_x
\end{displaymath}
where $x$ is regarded as a point both in $X$ and $X_1=X_{D_{X/Y}}$,
and $j: X_1 \to X$ is the second projection of the base change of
$X/Y$ over $Y_1= D_{X/Y}\to Y$. By (2) and (3) we can apply
\Proposition{\ref{spec-der-lemma}} to conclude that the map
\begin{displaymath}
  \psi : j^*(T_{X/Y}) \to T_{X_1} 
\end{displaymath}
is an isomorphism. Since $X$ is Cohen-Macaulay and $\pi$ is flat, it
follows that $X_1$ is Cohen-Macaulay, so that in particular $T_{X_1} $
satisfies $(S_2)$ \Lem{\ref{djuplemma}(1)}; therefore $ j^*(T_{X/Y})$
satisfies $(S_2)$. Since $\pi$ is flat, the dimension equality states that
\begin{displaymath}
\hto_{X_1}(x) + \hto_Y(y) =   \hto_X(x)\geq 3  
\end{displaymath}
when $x\in B_{X/Y}$ by (3). Then \thetag{\&\&} implies that
\begin{displaymath}
  \depth_{x}(T_{X/Y,x})\geq 3 \quad \text{when }\quad  x\in B_{X/Y}.
\end{displaymath}

It remains to prove the assertion in (5), where we can assume that
$Y= \Spec A$ and $X= \Spec B$, and $\pi: A\to B$ is a faithfully flat
morphism of local noetherian rings, where $B$, and therefore $A$ is
Cohen-Macaulay \cite{matsumura}*{\S 23 Corollary}. That $T_{B/A}$ is
flat when $\dim A \leq 1$ is evident, so assume that we have a regular
system of parameters $y_1, y_2$ for $A$. Since $B/A$ is faithfully
flat, $y_1,y_2$ is also a $B$ -sequence, and thus a $T_{B/A}$-sequence
\Lem{\ref{djuplemma}}. Since $A/Ay_2$ is a principal ring and
$T_{B/A}/y_2 T_{B/A}$ is a torsion free module, it is also flat over
$A/Ay_2$. Now apply \cite{matsumura}*{Exercise 22.3} to conclude that
$T_{B/A}$ is flat over $A$.
\end{proof}

\begin{pfof}{ the \hypertarget{assertion}{assertion} in
    \Remark{\ref{gensmooth}}}
  Since $X/Y$ is dominant and $X$ is regular at all points in the
  generic fibres, and the problem is local at such fibres, it follows
  that we can assume that $X$ and $Y$ are integral.  We will prove
  that if $\eta \in \Max (Y)$, then $\eta \not\in D_{X/Y}$.  Since
  $X/Y$ is dominant, there exists $\xi \in \Max(X)$ such that
  $\pi(x)=\pi(\xi)= \eta$, and we can moreover let $\xi$ be any
  maximal point that specialises to $x$, since it will satisfy
  $\pi(\xi)= \eta$ because $\eta$ is maximal.  We then have (as
  detailed below)
\begin{align*}
  \dim_{k_{X,x}}k_{X,x} \otimes_{\Oc_X,x}\Omega_{X/Y,x} & =
  \dim_{k_{X,x}} \mf_{X,x}/\mf_{X,x}^2+
  \dim_{k_{X,x}}\Omega_{k_{X,x}/k_{Y,\eta}} \\
  &= \hto_X(x) + \trdeg k_{X,x} /k_{Y,\eta} \\
  & = \hto_Y(\eta) + \trdeg k_{X,\xi}/k_{Y,\eta} \\
  & = \trdeg k_{X,\xi}/k_{Y,\eta} = \dim _{k_{X,\xi}}\Omega_{X/Y, \xi}
  = \rank \Omega_{X/Y,\xi}.
\end{align*}
The first line follows since $k_{X,x}/k_{Y,\eta}$ is separable, hence
0-smooth, after applying the second fundamental exact sequence in
\cite[Th. 25.2]{matsumura}. The second line follows since
$\Oc_{X,x}$ is regular and since $k_{X,x}/k_{Y,\eta}$ is separable and
finitely generated, so that a differential basis is the same as a
transcendence basis (see \cite[\S 26]{matsumura}).  The third line
follows since $X$ and $Y$ are integral, $X/Y$ is locally of finite
type, and $Y$ is N{\oe}therian, so that Ratliff's dimension equality holds
\cite[Th.  15.6]{matsumura}.

The see the second to last equality it suffices as above to prove that
$ k_{X,\xi}/k_{Y,\eta}$ is finitely generated and separable, where the
finite generation follows since $\pi$ is finitely presented at $x$.
First we note that $\Oc_{X,\xi} = k_{X,\xi}$ since it is a
localisation of the regular ring $\Oc_{X,x}$ and $\xi\in \Max (X) $.
Secondly, $k_{X,x}/k_{Y,\eta}$ is separable and $\Oc_{X,x}$ is
regular, hence $\Oc_{X,x}/k_{Y,\eta}$ is $\mf_{X,x}$-smooth
\cite[Lemma 1]{matsumura}.  Therefore the localisation
$\Oc_{X,\xi}/k_{Y,\eta}$ is $\mf_{X,\xi}$-smooth, hence $k_{X,\xi} $
is formally smooth over $k_{Y,\eta}$, implying that
$k_{X,\xi}/k_{Y,\eta}$ is separable \cite[Theorem 26.9]{matsumura}.

The last equality follows since $\Oc_{X,\xi}$ is regular, so
$\Oc_{X,\xi}= k_{X,\xi}$.  Since the first and the last entries are
equal it follows that $\Omega_{X/Y,x}$ is free (see ($*$)), so $\eta=
\pi(x) \not\in D_{X/Y}$.
\end{pfof}



\bigskip \textsc{Department of Mathematics, University of G\"avle,
  S-801 76, Sweden}

{\it E-mail address}\/: \texttt{rkm@hig.se}.

              \end{document}
              The following corollary is basic to the proof of \Theorem{\ref{Z-L}}. 

\begin{corollary}\label{can-iso}Assume that $Y$ is a regular scheme
  such that $1 \leq \dim Y \leq 2$, $X$ is connected and
  Cohen-Macaulay, and the relative dimension $d_{X/Y}\geq 2$.
  \begin{enumerate}
  \item Let $x$ be a point in $\Max (B_{X_1/Y_1})$ such that
    $\depth(x)\geq 2$, and $\pi$ is flat at $x$. Then we have an
    isomorphism
\begin{displaymath}
  \psi_x : j^*(T_{X/Y})_x \to T_{X_1/Y_1, x}. 
\end{displaymath}
\item Assume that $X/Y$ is flat and $\depth_{B_{X_1/Y_1}}(\Oc_{X_1})
  \geq 2$.  Then $\psi$ defines an isomorphism $ j^*(T_{X/Y}) \cong T_{X_1/Y_1}$.
  \end{enumerate}
\end{corollary}

\begin{lemma}\label{depth3} Let $A\to B$ be flat homomorphism of local
  rings, and put $\bar B = k_A\otimes_A B$.
  \begin{enumerate}
  \item Assume that $\depth B > \depth A$ and that $A$ is a
    domain. Then any $\bar B$-regular element is $k_A\otimes_A
    T_{B/A}$-regular, and $\depth k_A\otimes_A T_{B/A} \geq 1$.
  \item Assume that $A$ is regular, $\dim A\leq 2$, and $\depth B \geq
    2$. Then $T_{B/A}$ is flat over $A$.
  \item Assume that
    \begin{enumerate}
    \item $T_{B/A}$ is flat over $A$.
    \item $\depth A \geq 1$, $\depth B \geq  \depth A +2$, and $A$ is an
      integral domain.
    \end{enumerate}
    Then $\depth T_{B/A} \geq 3$, and $\depth k_A\otimes_A T_{B/A}\geq
    2$.
  \end{enumerate}
      \end{lemma}
      \begin{pf} (1): Since $A\to B$ is flat, $\depth k_A\otimes_A B =
        \depth B - \depth A \geq 1 $, so there exists a $k_A\otimes_A
        B$- regular element $\bar x \in k_A \otimes_A \mf_B $.  Take
        $\delta \in T_{B/A}$ and put $\bar \delta = 1\otimes \delta
        \in k_A\otimes_A T_{B/A}$, and assume that $\bar x \bar \delta
        =0$, so $x \delta = y \eta$, where $y\in \mf_A$, $x\in \mf_B$
        maps to $\bar x$, and $\delta \in T_{B/A}$. Since $A \to B$ is
        flat and $A$ is a domain, it follows that $y,x$ is a
        $B$-sequence (see the proof of \cite{matsumura}*{Th. 23.3}),
        hence it is a $T_{B/A}$-sequence \Lem{\ref{djuplemma}};
        therefore $\delta = y \delta'$ for some $\delta'\in T_{B/A}$,
        so $\bar \delta =0$, and hence $\bar x$ is $k_A\otimes_A
        T_{B/A}$-regular.

        (2): The assertion is evident when $\dim A \leq 1$, so assume
        that we have a regular system of parameters $y_1, y_2$ for
        $A$.  By faithful flatness it is also a $B$ -sequence, and
        thus a $T_{B/A}$-sequence \Lem{\ref{djuplemma}}.  Since
        $A/Ay_2$ is a principal ring and $T_{B/A}/y_2 T_{B/A}$ is a
        torsion free module, it is also flat over $A/Ay_2$. Now apply
        \cite{matsumura}*{Exercise 22.3} to conclude that $T_{B/A}$ is
        flat over $A$.

        (3): When $\depth A \geq 2$, \cite{matsumura}*{Th. 23.3} and
        (1) implies that $\depth T_{B/A}\geq 3$. If $\depth A =1$,
        since $A$ is normal, $A$ is smooth, so that $\mf_A $ is a
        principal ideal, hence we are in the situation of
        \Proposition{\ref{scheja-storch}}. Thus since $T_{\bar B/k}$
        is free, it follows that $\bar B/k$ is smooth; since $B/A$ is
        flat, it follows that $B/A$ is smooth; therefore $T_{B/A}$ is
        locally free, so that $\depth T_{B/A} \geq 3$ if and only of
        $\depth B \geq 3$. By [loc cit, Corollary], there exist
        $\bar x_1, \bar x_2\in \mf_{\bar B}$ forming a $\bar B$-
        sequence, and since $A$ is an integral domain, any element
        $y\in \mf_A$ is $A$-regular. Let $x_1, x_2\in \mf_B$ be
        elements that map to $\bar x_1 $ and $\bar x_2$, respectively.
        By the argument in the proof of \cite{matsumura}*{Th. 23.3 },
        it follows that $x_1, x_2, y$ forms a $B$-sequence.
      \end{pf}

\begin{pfof}{\Proposition{\ref{can-iso}}}
  (1): Since the assertion is about the stalk at a point $x$ we can
  assume that $X= \Spec B$ and $Y = \Spec A$, and $\pi$ is defined by
  a flat homomorphism of local rings $A \to B$. The base change in the
  diagram (BC) arises from a homomorphism of local rings $A \to A_1$,
  so $Y_1= \Spec A_1$; let $x$ be the closed point in $X_1$. Put
  $X^0 = X \setminus \{j (x)\}$ and $X_1^0 = X_1 \setminus \{ x \} $.
  Let $\phi_1: X_1^0 \to X_1$ and $ j_0 : X^0_1 \to X^0 $ be the
  canonical maps occuring in the base change of $\phi: X^0 \to X $
  over $j$, so that $\phi\circ j_0 = j\circ \phi_1$. Now in general,
  for a coherent $\Oc_X$-module we have a base change homomorphism
  $j^*\phi_*(M)\to (\phi_1)_* j_0^*(M)$ (sections in punctured
  neighbourhoods of $x$ in $X$ are mapped to sections in a punctured
  neighbourhood of $x$ in $X_1$ ), which need neither be injective nor
  surjective. However, we have:
  \begin{enumerate}[label=({\roman*})]
  \item $\Omega_{X_0/Y,x'}$ is free when $x' \in j(X^0)$, so
    $ \phi^*(T_{X/Y})_{x'} = T_{X^0_1/Y_1, x'}$ when $x'\in X^0$, and
    $j_0^*\phi^*(T_{X/Y}) = \phi_1^*(T_{X_1/Y_1})$.
  \item Since $\hto_{X}(j(x))\geq 2$ and $\Oc_{X,x}$ is torsion free,
    $ \phi_*\phi^*(T_{X/Y})_x = T_{X/Y,x}$ \Lem{\ref{djuplemma}}.

    Since $Y$
    is regular of dimension $\leq 2$, \Lemma{\ref{depth3}}, (2),
    implies that $T_{X/Y,j(x)}$ is flat over $\Oc_{Y,\pi(j(x))}$,
    hence by \Lemma{\ref{depth3}, (3)}, $\depth j^*(T_{X/Y})_x\geq 2$;
    therefore ${i_1}_*i_1^*j^* ( i_*(T_{X/Y}))_x = j^* (
    i_*(T_{X/Y}))_x =0$. Since $\depth (x)\geq 2$, so $\depth
    T_{X_1/Y_1}\geq 2$ \Lem{\ref{djuplemma}}, we also have
    $(i_1)_*i_1^*(T_{X_1/Y_1})_x = T_{X/Y,x}$. 
\end{enumerate} 
We then get from (i-ii), and since $(i_1)^*j^* = j_0^*i^*$,
    \begin{eqnarray*}
      j^*  (  T_{X/Y})_x &=& (i_1)_*i_1^* (j^* (
      T_{X/Y}))_x=  (i_1)_* ( j_0^*i^*(T_{X/Y}))_x \\ &=& (i_1)_*
      i_1^*(T_{X_1/Y_1})_x  = T_{X_1/Y_1,x}.
  \end{eqnarray*}

  (2): Clearly, $\psi_x$ is an isomorphism when $\hto_{X_1}(x)\leq 1$,
  and by (1) it is also an isomorphism when $\hto_{X_1}(x)\geq 2$.
\end{pfof}

Although it is not used elsewhere, we include the following purity
result.

Let $T_{X/S\to Y/S}$ be the sheaf of derivations $\pi^{-1}(\Oc_Y)\to
\Oc_X$, so we have a tangent morphism of $\Oc_X$-modules $d\pi :
T_{X/S}\to T_{X/S\to Y/S}$. Similarly, $T_{A/k \to B/k}= \Gamma(\Spec
B, T_{\Spec B/\Spec k \to \Spec A/\Spec k})$ denotes the $B$-module of
$k$-linear derivations $A\to B$.

\begin{prop}\label{lemmasurjective}
  Let $(A, \mf_A)\to (B, \mf_B)$ be a flat homomorphism of local
  rings, $k $ is a subring of $A$, and $\pi : X= \Spec B \to Y= \Spec
  A$ its associated map of affine $S$-schemes, where $S= \Spec k$.
  Assume:
\begin{enumerate}
\item  The discriminant set $\pi(B_{X/Y}) = \{y\}$, the closed point
  in $Y$.
\item $\codim^-_{X_y} B_{X_y/k_{Y,y}}\geq 1$.
\item  $\depth B \geq 3$, $\dim A =2$, and $A$ is
  regular.
\end{enumerate}
Then the tangent map $d \pi : T_{B/k}\to T_{A/k \to B/k}$ is
surjective.
\end{prop}
\begin{proof}
  Put $y=\pi(x)$, let $X_{y}\subset X$ be the special fibre, put
  $X_{0}= X \setminus B_{X/Y}$, and let $i: X_{0}\to X$ be the open
  inclusion.  $\pi$ is flat and $\dim A \geq 2$, so $\codim^-_X
  {B_{X/Y}}\geq 2$, and since $B$ satisfies $(S_2)$, $\depth_{B_{X/Y}}
  T_{X/S} \geq 2$ and $B$ is torsion free so $ T_{X/S}$ is torsion
  free. Therefore $i_*i^{*}(T_{X/S}) = T_{X/S}$. Similarly, $T_{X/S\to
    Y/S } = i_*(T_{X_{0}/S\to Y/S})$.

  We have an exact sequence of sheaves of $\Oc_{X_{0}}$-modules
    \begin{displaymath}
0 \to T_{X_{0}/Y} \to T_{X_{0}/S}\to T_{X_{0}/S   
      \to Y/S}\to 0
  \end{displaymath}
  and hence the exact sequence
\begin{multline*}
  0 \to \Gamma(X, i_*i^*(T_{X/Y}))\to \Gamma(X, i_*i^*(T_{X/S}\to T_{X/S   
      \to Y/S})) \\ \to H^1(X_{0}, T_{X/Y}) \to   
  \end{multline*}
  By the previous paragraph one can erase $i_*i^*$ from this sequence.
  If $x\in B_{X/Y}$, so by (2) $\hto_{X_y}(x)\geq 1$, the dimension
  equality for flat morphisms implies that $\hto_X(x) = \hto_{X_y}(x)
  + \hto_Y(y)\geq 3$.  Then by the first line in the proof of
  \Lemma{\ref{depth3}}, (3), $\depth_{B_{X/Y}} (T_{X/Y})\geq 3$, and
  therefore $H^1(X_{0}, T_{X/Y})=0$.  This implies that the tangent
  map $T_{B/k}\to T_{A/k \to B/k}$ is surjective.
\end{proof}
We have a kind of converse of \Proposition{\ref{lemmasurjective}}.
\begin{prop}\label{depth**}Assume:
  \begin{enumerate}
  \item The tangent morphism $d\pi : T_{B/k}\to T_{A/k \to B/k}$ is
    surjective.
  \item $\depth B \geq 3$ and any $B$-sequence of length 3 is also a $
    T_{B/k}$-sequence (e.g. $B/k$ is smooth)
\item $\depth k_A\otimes_A B \geq 2$ (e.g $A\to B$ is flat and $\depth
  B \geq 2 + \depth A$)
\end{enumerate}
Then:
\begin{enumerate}
\item If $x_1,x_2, x_3\in \mf_B$ is a $T_{B/k}$-sequence, then it is
  also a $T_{B/A}$- sequence.
\item If $A\to B$ is flat and $A$ is an integral domain, then $\depth
  k_A \otimes_A T_{B/A}\geq 2$.
\end{enumerate}
\end{prop}
\begin{proof}
  (1): If $x_1, x_2\in \mf_B$ is $B$-sequence, then it is also a
  $T_{B/A}$ - and $T_{B/k}$-regular sequence \Lem{\ref{djuplemma}},
  and since $\depth B \geq 3$, we can find an element $x_3\in\mf_B$
  such that $x_1, x_2, x_3\in \mf_B$ forms a
  $T_{B/k}$-sequence. Hence, if $\delta_1, \delta_2, \delta_3 \in
  T_{B/A}$ satisfies $x_1 \delta_1 + x_2 \delta_2 + x_3\delta_3=0$,
  there exist $\partial_1, \partial_2 \in T_{B/k}$ such that $\delta
  _3 = x_1 \partial_1 + x_2 \partial_2$.  If $a\in A $ we then have
  \begin{displaymath}
    x_1 \partial_1(a) +  x_2 \partial_2(a) =0,
  \end{displaymath}
  and since $x_1, x_2$ is a $B$-sequence, we have $\partial_1(a) = x_2
  b $ and $\partial_2(a)= -x_1 b$, for some element $b\in B$. Since
  $x_1$ and $x_2$ are $B$-regular, it follows that the map $\eta:
  a\mapsto b $ defines a derivation $ A\to B$, so $\eta \in T_{A/k \to
    B/k}$.  Since $d\pi$ is surjective there exists an element $\hat
  \eta \in T_{B/k}$ such that $d\pi (\hat \eta) = \eta $. Then putting
  $\delta'_1 = \partial_1 - x_2 \hat \eta$ and $\delta'_2 = \partial_2
  + x_1 \hat \eta$, we have $\delta_1', \delta_2' \in T_{B/A}$ and
  $\delta_3 = x_1 \delta_1' + x_2 \delta_2'$.

  (2): Since $T_{B/A}$ is flat over $A$, the assertion follows from
  (1) and \cite{matsumura}*{Th 23.3}.
\end{proof}
\begin{lemma}\label{2dimlemma} Let $\pi : X\to Y$ be a morphism as in
  \Theorem{\ref{Z-L}}. Given a point $y\in Y$ of height $\hto_Y(y)\geq
  1$, there exists a (restriction of a) base change $\pi': X' \to Y' $
  which satisfies the conditions in \Proposition{\ref{can-iso}} and
  $X' $ is Cohen-Macaulay, and there is a point $y'$ that maps to $y$
  such that the fibre $X'_{y'}$ is a base change of $X_y$ over a
  finite field extension.
\end{lemma}
\begin{pf}
  If $\hto_Y(y) = 1$, we let $Y' = \Spec \Oc_{Y,y}$ and $Y' \to Y$ be
  induced by the restriction morphism, and $y'=y$. The base change
  $X'\to Y'$ then remains a flat d.c.i.  and the special fibre $X'_y$
  coincides with $X_y$. Now assume $\hto_Y(y)\geq 2$, and select a
  point $y_1$ that specialises to $y $ and $\hto_Y(y') = \hto_Y(y)
  +1$. Let $Y_1$ be the normalisation of the closure ${y_1}^-$, so we
  have a finite surjective morphism $Y_1 \to Y$.  The base change
  $X_1\to Y_1$ remains a flat d.c.i., and the fibre $(X_1)_{y_1}$ over
  a point $y_1$ that maps to $y$ is a base change of $X_y$ over a
  finite field extension. Let now $Y' = \Spec \Oc_{Y_1, y_1}$ and $Y'
  \to Y_1$ be induced by restriction of functions, and $X'\to Y'$ be
  the corresponding base change of $X_1\to Y_1$.

  In either case $\hto_Y(y)\geq 1$, clearly $X'/Y'$ is flat and also a
  d.c.i., hence $X'_{y'}$ is a local complete intersections and
  therefore Cohen-Macaulay; hence $X'/Y'$ is Cohen-Macaulay at each
  point in $X'_y \subset X'$, and since the set of points where $X'$
  is Cohen-Macaulay is open, we can restrict, if necessary, the
  morphism $X' \to Y'$ to such an open subset of $X'$, containing
  $X_{y'}' $ , so that $X'$ can be assumed to be Cohen-Macaulay.
\end{pf}

By \Proposition{\ref{lipman-normal}} it suffices to prove $\codim_X^+
B_{X/k}\leq 1$.  Let $j:X/k\to Z/k$ be a regular immersion into a
smooth variety $Z/k$, so locally $j(X)$ is defined by an ideal $I$
such that $I/I^2$ is locally free over $\Oc_X$ and we have the short
exact sequence $0\to I/I^2 \to j^*(\Omega_{Z/k})\to \Omega_{X/k}\to
0$, and in particular $\pdo \Omega_{X/k,x}\leq 1$ for each point $x$
in $X$. Dualising we get the exact sequences
\begin{align*}
  0&\to T_{X/k} \to j^*(T_{Z/k})\to \Cc_{X/Z}\to 0, \tag{E} \\
 0 &\to \Cc_{X/Z} \to (I/I^2)^* \to Ext_{\Oc_X}^1(\Omega_{X/k},\Oc_X)\to 0.\tag{F}
\end{align*}

(a) The normal module $\Cc_{X/Z}$ satisfies $\depth \Cc_{X/Z,x}\geq 2$
when $I$ is a locally principal ideal, $\hto_X(x)\geq 2$, and
$x\in \Max (B_{X/k})$:

A locally defined surjection
$\Oc_Z ^{d_X}\to j_*(T_{X/k})$ together with the surjective map
$T_{Z/k}\to j_{*}j^*(T_{Z/k})$ gives a lift
\begin{displaymath}\tag{E'}
0 \to \Oc_Z^{d_X}\to T_{Z/k} \to \hat \Cc_{X/Z}\to 0
\end{displaymath}
of the sequence \thetag{$E$}, which is exact to the left since
$T_{X/k}$ is locally free. Since the ideal $I\subset \Oc_{Z}$ is
locally principal, it follows that the cokernel $\hat \Cc_{X/Z}$ is an
ideal of $\Oc_X$, and
$\hat \Cc_{X/Z,j(x)} = T_{Z/k,j(x)}\cdot f \subset \Oc_{Z,j(x)}$, when
$f$ is a local generator which is defined near $j(x)$. We note also
that the element $f$ belongs to the integral closure of the ideal
$T_{Z/k,j(x)}\cdot f$; this is clear when $\hto_Z(j(x)) \leq 1$ and
the general case follows from describing the integral closure of an
ideal $J_{z} \subset \Oc_{Z,z}$ as the intersection
$\cap J_z \Oc_{Z,z'}$, running over points $z'$ such that
$\hto_Z(z')\leq 1$ and $z'$ specialises to $z$. Since
$\pdo \hat \Cc_{Z/X,j(x)} \leq 1$, by the Hilbert-Burch theorem
$ \hat \Cc_{X/Z ,j(x)} = a F_1(\hat \Cc_{X/Z,j(x)})$ for some
$a\in \Oc_{Z,j(x)}$; since $x\in \Max(B_{X/k})$ and $\hto_X(x)\geq 2$;
hence by Krull's principal ideal theorem $a$ restricts to a unit in
$\Oc_{X,x}$ and therefore it is a unit in $\Oc_{Z,j(x)}$; therefore
\begin{displaymath}
  T_{Z/k,j(x)}\cdot f = \hat \Cc_{Z/X,j(x)}= F_1(\hat
  \Cc_{X/Z,j(x)}).
\end{displaymath}
( {\underline Note}: The weaker assertion
$V_Z(T_{Z/k,j(x)}\cdot f)=V(F_1(\hat \Cc_{X/Z,j(x)}))$ is actually
sufficient for the proof, and is easy to see: the germ $ V(\hat
\Cc_{X/Z,j(x)}) $ is of codimension $\geq 2$ and therefore $\hat
\Cc_{X/Z,z}$ is not principal if and only if $z\in V(\hat
\Cc_{X/Z,j(x)})$, by Krull's principal ideal theorem.)  Considering
germs of varieties at $x$ we now get
\begin{align*}
V_X(F_1(\Cc_{X/Z})) & =V_{Z}(F_1(\hat \Cc_{X/Z}) +I) = V_Z(T_{Z/k}\cdot f +I)\\ & = V_Z(T_{Z/k}\cdot
f) = V_{Z}(F_1(\hat \Cc_{X/Z})),
\end{align*}
so in particular $\Max (V_X(F_1(\Cc_{X/Z})))= \Max ( V_Z(F_1(\hat
\Cc_{X/Z})))$, identifying $X$ with $j(X)$.  By the Eagon-Northcott
bound on heights of determinant ideals \cite{eagon-northcott:height}
applied to the exact sequence \thetag{E'}, we get $\hto_Z(z)\leq 2$
when $z\in \Max (V_Z(F_1(\hat \Cc_{Z/X})))$ and $z$ specialises to
$j(x)$.  Therefore, if $j(x)\in \Max( V_X(F_1(\Cc_{X/Z})))$, we get
$\hto_X(x) = \hto_Z(j(x)) -1 \leq 1$.  Since $\hto_X(x) \geq 2$ it
follows that $\Cc_{X/Z,x}$ is free, and since $X$ is Cohen-Macaulay,
$\depth \Cc_{X/Z,x} \geq 2$.

(b) $\codim^+_X B_{X/k}\leq 1$: Assume now on the contrary that there
exists a point $x\in \Max (B_{X/k})$ such that $\hto (x)\geq 2$. By
(a)
\begin{displaymath}
Ext^1_{\Oc_X,x}(Ext^1_{\Oc_{X,x}}(\Omega_{X/k,x},\Oc_{X,x}),\Cc_{X/Z,x}) =0,
\end{displaymath}
hence the sequence \thetag{F} splits, so there exists an injection
$Ext^1_{\Oc_{X,x}}(\Omega_{X/k,x}, \Oc_{X,x})$ $ \to (I_x/I_x^2)^*$.
Since $X$ is Cohen-Macaulay, the free module $(I_x/I_x^2)^*$ has no
embedded associated prime. Therefore
\begin{displaymath}
  Ext^1_{\Oc_{X,x}}(\Omega_{X/k,x}, \Oc_{X,x}) =0.
\end{displaymath}
Since $\pdo \Omega_{X/k,x}\leq 1$, this implies that $\Omega_{X/k,x}$
is free, contradicting the assumption that $x\in B_{X/k}$. Therefore
$\codim^+_X B_{X/k}\leq 1$.

\begin{pf} (1): Since $A\to B$ is flat,
  $\depth k_A\otimes_A B = \depth B - \depth A \geq 1 $, so there
  exists a $k_A\otimes_A B$- regular element
  $\bar x \in k_A \otimes_A \mf_B $. Take $\delta \in T_{B/A}$ and put
  $\bar \delta = 1\otimes \delta \in k_A\otimes_A T_{B/A}$, and assume
  that $\bar x \bar \delta =0$, so $x \delta = y \eta$, where
  $y\in \mf_A$, $x\in \mf_B$ maps to $\bar x$, and
  $\delta, \eta \in T_{B/A}$. Since $A \to B$ is flat and $A$ is a
  domain, it follows that $y,x$ is a $B$-sequence (see the proof of
  \cite{matsumura}*{Th. 23.3}), hence it is a $T_{B/A}$-sequence
  \Lem{\ref{djuplemma}}; therefore $\delta = y \delta'$ for some
  $\delta'\in T_{B/A}$, so $\bar \delta =0$, and hence $\bar x$ is
  $k_A\otimes_A T_{B/A}$-regular.

        (2): The assertion is evident when $\dim A \leq 1$, so assume
        that we have a regular system of parameters $y_1, y_2$ for
        $A$.  By faithful flatness it is also a $B$ -sequence, and
        thus a $T_{B/A}$-sequence \Lem{\ref{djuplemma}}.  Since
        $A/Ay_2$ is a principal ring and $T_{B/A}/y_2 T_{B/A}$ is a
        torsion free module, it is also flat over $A/Ay_2$. Now apply
        \cite{matsumura}*{Exercise 22.3} to conclude that $T_{B/A}$ is
        flat over $A$.

        (3): When $\depth A \geq 2$, \cite{matsumura}*{Th. 23.3} and
        (1) implies that $\depth T_{B/A}\geq 3$. If $\depth A =1$,
        since $A$ is normal, $A$ is smooth, so that $\mf_A $ is a
        principal ideal, hence we are in the situation of
        \Proposition{\ref{scheja-storch}}. Thus since $T_{\bar B/k}$
        is free, it follows that $\bar B/k$ is smooth; since $B/A$ is
        flat, it follows that $B/A$ is smooth; therefore $T_{B/A}$ is
        locally free, so that $\depth T_{B/A} \geq 3$ if and only of
        $\depth B \geq 3$. By [loc cit,\S 23 Corollary], there exist
        $\bar x_1, \bar x_2\in \mf_{\bar B}$ forming a $\bar B$-
        sequence, and since $A$ is an integral domain, any element
        $y\in \mf_A$ is $A$-regular. Let $x_1, x_2\in \mf_B$ be
        elements that map to $\bar x_1 $ and $\bar x_2$, respectively.
        By the argument in the proof of \cite{matsumura}*{Th. 23.3 },
        it follows that $x_1, x_2, y$ forms a $B$-sequence.
      \end{pf}
\begin{proof}
  Put $y=\pi(x)$, let $X_{y}\subset X$ be the special fibre, put
  $X_{0}= X \setminus B_{X/Y}$, and let $i: X_{0}\to X$ be the open
  inclusion.  $\pi$ is flat and $\dim A \geq 2$, so $\codim^-_X
  {B_{X/Y}}\geq 2$, and since $B$ satisfies $(S_2)$, $\depth_{B_{X/Y}}
  T_{X/S} \geq 2$ and $B$ is torsion free so $ T_{X/S}$ is torsion
  free. Therefore $i_*i^{*}(T_{X/S}) = T_{X/S}$. Similarly, $T_{X/S\to
    Y/S } = i_*(T_{X_{0}/S\to Y/S})$.

  We have an exact sequence of sheaves of $\Oc_{X_{0}}$-modules
    \begin{displaymath}
0 \to T_{X_{0}/Y} \to T_{X_{0}/S}\to T_{X_{0}/S   
      \to Y/S}\to 0
  \end{displaymath}
  and hence the exact sequence
\begin{multline*}
  0 \to \Gamma(X, i_*i^*(T_{X/Y}))\to \Gamma(X, i_*i^*(T_{X/S}\to T_{X/S   
      \to Y/S})) \\ \to H^1(X_{0}, T_{X/Y}) \to   
  \end{multline*}
  By the previous paragraph one can erase $i_*i^*$ from this sequence.
  If $x\in B_{X/Y}$, so by (2) $\hto_{X_y}(x)\geq 1$, the dimension
  equality for flat morphisms implies that
  $\hto_X(x) = \hto_{X_y}(x) + \hto_Y(y)\geq 3$. By
  \Lemma{\ref{depth3}}(3), $\depth_{B_{X/Y}} (T_{X/Y})\geq 3$, and
  therefore $H^1(X_{0}, T_{X/Y})=0$. This implies that the tangent map
  $T_{B/k}\to T_{A/k \to B/k}$ is surjective.
\end{proof}
Tthere is also a kind of converse of
\Proposition{\ref{lemmasurjective}}.
\begin{prop}\label{depth**}Assume:
  \begin{enumerate}
  \item The tangent morphism $d\pi : T_{B/k}\to T_{A/k \to B/k}$ is
    surjective.
  \item $\depth B \geq 3$ and any $B$-sequence of length 3 is also a $
    T_{B/k}$-sequence (e.g. $B/k$ is smooth)
\item $\depth k_A\otimes_A B \geq 2$ (e.g $A\to B$ is flat and $\depth
  B \geq 2 + \depth A$)
\end{enumerate}
Then:
\begin{enumerate}
\item If $x_1,x_2, x_3\in \mf_B$ is a $T_{B/k}$-sequence, then it is
  also a $T_{B/A}$- sequence.
\item If $A\to B$ is flat and $A$ is an integral domain, then $\depth
  k_A \otimes_A T_{B/A}\geq 1$.
\end{enumerate}
\end{prop}
\begin{proof}
  (1): If $x_1, x_2\in \mf_B$ is $B$-sequence, then it is also a
  $T_{B/A}$ - and $T_{B/k}$-regular sequence \Lem{\ref{djuplemma}},
  and since $\depth B \geq 3$, by (2) we can find an element
  $x_3\in\mf_B$ such that $x_1, x_2, x_3\in \mf_B$ forms a
  $T_{B/k}$-sequence. Hence, if
  $\delta_1, \delta_2, \delta_3 \in T_{B/A}$ satisfies
  $x_1 \delta_1 + x_2 \delta_2 + x_3\delta_3=0$, there exist
  $\partial_1, \partial_2 \in T_{B/k}$ such that
  $\delta _3 = x_1 \partial_1 + x_2 \partial_2$. If $a\in A $ we then
  have
  \begin{displaymath}
    x_1 \partial_1(a) +  x_2 \partial_2(a) =0,
  \end{displaymath}
  and since $x_1, x_2$ is a $B$-sequence, we have $\partial_1(a) = x_2
  b $ and $\partial_2(a)= -x_1 b$, for some element $b\in B$. Since
  $x_1$ and $x_2$ are $B$-regular, it follows that the map $\eta:
  a\mapsto b $ defines a derivation $ A\to B$, so $\eta \in T_{A/k \to
    B/k}$.  Since $d\pi$ is surjective there exists an element $\hat
  \eta \in T_{B/k}$ such that $d\pi (\hat \eta) = \eta $. Then putting
  $\delta'_1 = \partial_1 - x_2 \hat \eta$ and $\delta'_2 = \partial_2
  + x_1 \hat \eta$, we have $\delta_1', \delta_2' \in T_{B/A}$ and
  $\delta_3 = x_1 \delta_1' + x_2 \delta_2'$.

  (2): Since $T_{B/A}$ is flat over $A$ \Lem{\ref{depth3}(2)}, the assertion follows from
  (1) and \cite{matsumura}*{Th 23.3}.
\end{proof}


